 \documentclass[11pt,letterpaper]{amsart}
\usepackage{graphicx}
\usepackage[all]{xy}
\usepackage{amsmath}
\usepackage{amsthm}
\usepackage{latexsym,bm}
\usepackage{amssymb}
\usepackage{indentfirst}
\newtheorem{thm}{Theorem}[section]
\newtheorem{cor}[thm]{Corollary}
\newtheorem{lem}[thm]{Lemma}
\newtheorem{prop}[thm]{Proposition}
\newtheorem{defn}[thm]{Definition}
\newtheorem{rem}[thm]{Remark}
\newtheorem{exm}{Example}
\numberwithin{equation}{section}

\begin{document}

\title[New Laplacian comparison \&. Finsler-grad-estimates]{\bf{\Large Operators on nonlinear metric measure spaces I: A New Laplacian comparison theorem on Finsler manifolds and a generalized Li-Yau approach to gradient estimates of Finslerian Schr\"odinger equation}}

\author{{\bf Bin Shen}}

\date{}

\maketitle

\begin{quote}
\small {\bf Abstract}. This paper is the first in a series of paper where we describe the differential operators on general nonlinear metric measure spaces, namely, the Finsler spaces. We try to propose a general method for gradient estimates of the positive solutions to nonlinear parabolic (or elliptic) equations on both compact and noncompact forward complete Finsler manifolds. The key of this approach is to find a valid comparison theorem, and particularly, to find related suitable curvature conditions. Acturally, we define some non-Riemannian tensors and a generalization of the weighted Ricci curvature, called the mixed weighted Ricci curvature. With the assistance of these concepts, we prove a new kind of Laplacian comparison theorem. 
As an illustration of such method, we apply it to Finslerian Schr\"odinger equation.
The Schr\"odinger equation is a typical class of parabolic equations, and the gradient estimation of their positive solution is representative. In this paper, based on the new Laplacian comparison theorem, we give global gradient estimates for positive solutions to the Finslerian Schr\"odinger equation in $CD(K,N)$ conditions, and local gradient estimates with respect to two different mixed weighted Ricci curvature conditions. More interestingly, we found that the $CD(K,N)$ condition was used in the generalization of Li-Yau type estimations, and the $CD(K,\infty)$ condition was used in the processing of Hamilton estimates.

\end{quote}

\begin{quote}
\end{quote}

\baselineskip 17pt



\section{Introduction}

For nearly four decades, scholars have studied problems related to geometric analysis on different differential manifolds. As far as the topic of this article is concerned, the work on gradient estimation is also endless. In 1986, P. Li and S-T. Yau made important research on operator kernel estimation of Schr\"odinger equations on manifolds. They startied from gradient estimation, gave Hanack's inequality and the upper and lower bounds of the fundamental solution, and studied the spectral and isolatence properties of operators \cite{LiYau1986}. 
Along the method of their article, a great deal of gradient estimates have been given on Riemannian manifolds, as well as the Liouville type theorem and the blow up property of the solutions, etc. Some of the more recent work including \cite{Zhu2011NA,HHL2013AGAG,CKL2011JDG} and so on. Moreover, in terms of improvement of the curvature conditions, Qi.S. Zhang and M. Zhu adopted some $L^p$ integral Ricci curvature condition to get the local gradient eatimates \cite{ZZ2017Proceeding}. Later, they released the curvature condition allowing Ricci curvature not bounded form below \cite{ZZ2019JFA}. \cite{O2019proceeding} also used the same integral curvature to get the Neumann Li-Yau gradient estimates. C. Rose derived a Li-Yau gradient estimate under a Kato class condition for the negative part of the Ricci curvature over compact Riemannian manifolds \cite{R2019AGAG}.

In addition to gradient estimates of positive solutions to various equations on Riemannian manifolds, as a direct generalization, the weighted Riemannian space is a broader form of the Riemannian manifold to consider the geometric analysis problems. Actually, M. Erbar, K. Kuwada and K-T. Sturm proved the equivalence of the curvature-dimension bounds of Lott-Sturm-Villani and Bakry-\'Emery in complete generality for infinitesimally Hilbertian metric measure spaces \cite{EKS2015Invent}. R. Jiang built Moser-Trudinger and Sobolev inequalities of gradients for solutions to Poisson's equations \cite{Jiang2011JFA}. Local Li-Yau's estimate for weak solutions to the heat equation and a sharp Yau's gradient for harmonic functions on metric measure spaces were established by H. Zhang and X. Zhu in \cite{ZZ2016CV}. J. Wu established local elliptic and parabolic gradient estimates for positive smooth solutions to a parabolic equation with two more nonlinear factors than the Schr\"odinger equation \cite{Wu2018MM}. Local Hamilton-type gradient estimates and Harnack inequalities for positive solutions to another nonlinear equation more general than the Schr\"odinger equation were derived by F. Zeng in \cite{Zeng2021MJM}. H. Dung, N. Dung and J. Wu proved sharp gradient estimates for positive solutions to the weighted heat equation on smooth metric measure spaces with compact boundary \cite{DungDW2021CPAA}.  In \cite{AAM2022JGP}, A. Abolarinwa, A. Ali and F. Mofarreh finally generalized such results to a class of more general equations on complete metric measure spaces with or without boundary. Gradient estimates for a more nonlinear parabolic equation on metric measure spaces under geometric flow were provided by A. Shahroud \cite{Shah2023}.

The applications of gradient estimation on the K\"ahler manifolds have also formed an important subtopic in recent decades. L. Ni got an integral gradient estimate of the subharmonic function bounded from above and below on a complete noncompact K\"ahler manifold of complex dimension greater than 2 with bounded nonnegative holomorphic bisectional curvature, to prove some rigidity theorem \cite{Ni1998JDG}. Z. Lu, Q. Zhang and M. Zhu proved a certain gradient and eigenvalue estimates, as well as the heat kernel estimates, for the Hodge Laplacian on $(m, 0)$ forms, by using a new Bochner type formula for the gradient of $(m, 0)$ forms only involving the Ricci curvature and the gradient of the scalar curvature \cite{LZZ2021JGA}. 

Discrete Li-Yau gradient estimate for graphs satisfying the CDE-inequality was first presented by F. Bauer, P. Horn, Y. Lin, G. Lippner, D. Mangoubi and S-T. Yau \cite{BHLLMY2015JDG}. So far, the research of gradient estimation on graphs have focused mainly on the heat equation. L. Wang and Y. Zhou gave gradient estimates of the Laplacian eigenequations on locally finite infinite graphs \cite{WZ2017AM}. The gradient estimation of the positive solution to the heat equation on the finite graphs was determined by P. Horn \cite{H2019SIAM}. C. Gong, Y. Lin, S. Liu and S-T. Yau derived Li-Yau inequality for unbounded Laplacian on complete weighted graphs with the assumption of the curvature dimension inequality $CDE'(n,K)$. They also got Harnack inequality, heat kernel bounds and Cheng's eigenvalue estimate as the applications \cite{GLLY2019AM}. 

In the past decade, many scholars have also turned to the study of gradient estimations on differential manifolds whose metrics are evolved by a certain curvature flow. H. Cao and M. Zhu proved the Li-Yau-Hamilton and Aronson-B\'enilan estimates for fast diffusion and porous medium equations under the Ricci flow \cite{CZ2015JMPA,CZ2018NA}. A. Abolarinwa studied the global Harnack estimates for heat equations with potentials \cite{Abolar2016MJM}. L. Zhang investigated Li-Yau's type and J. Li's type parabolic gradient estimates for positive solutions to a nonlinear parabolic equation which is a nonlinear generalization of the Schr\"odinger equation under the Ricci flow \cite{Zhang2021BBMS}. He also researched the gradient estimates for another nonlinear PDE under the Yamabe flow \cite{Zhang2019JMAA}. There are more results on Riemannian manifolds under the geometric flow \cite{Sun2011PJM,DMT2023JMAA,LBZ2016JMAA,LZ2020PA}.\\

Finsler metric measure space could be considered as a nonlinear generalization of the weighted Riemannian space. The analysis on it is related to the optimal transport on irreversible metric spaces \cite{KZ2022JMPA} and uncertainty principles \cite{HKZ2020TAMS}, etc. Nowadays, research on gradient estimates on Finsler manifolds is still in its infancy. S. Ohta first proved the gradient estimate of heat equations on compact manifolds by establishing the curvature dimension condition on Finsler measure spaces \cite{Ohta2014}.
Based on their successful promotion of traditional research methods of the linear Laplacian, some scholars have obtained the analysis results on Finsler manifolds using curvature dimension conditions. In 2019, F. Zeng and Q. He studied the gradient estimation of the positive solution to the heat equation under the evolving of the metric along the Finsler-Ricci flow \cite{ZH2019MS}. Later, X. Cheng obtained a positive solution gradient estimate of the heat equation on a compact Finsler manifold under the Finsler-Ricci flow with more different curvature conditions \cite{Cheng2022JMAA}. S. Azami proved the gradient estimation of the positive solution to the heat equation under the general Finsler geometric flow \cite{Azami2022MS}. More generally, he subsequently provided gradient estimates and differential Harnack inequalities for positive solutions to a kind of nonlinear parabolic equation on compact Finsler manifolds, whose metric evolved by the general Finsler geometric flow \cite{Azami2023MJM}. All those results require that Finsler manifolds be compact. On a forward complete noncompact Finsler manifold, S. Ohta etc. \cite{Ohta2014} also asked the question: how to build a correlative theory on noncompact Finsler manifolds?  
C. Xia first adopted the Moser's iteration to investigate the positive harmonic function on noncompact Finsler manifolds \cite{XiaC2014CVPDE}. 
It is a good answer to Ohta's question in the harmonic function case, employing the same curvature conditions. 
Inspired by this work and using the Moser's iteration, Q. Xia obtained a Li-Yau type gradient estimation of the positive solution to the heat equation on forward complete noncompact Finsler manifold \cite{Xia2023}.  
However, the approach in \cite{Xia2023} has limitations 
(cf. Remark 4.1 in \cite{Xia2023}).

One might wonder why local gradient estimates can not be obtained on noncompact Finsler manifolds using the same method as on Riemannian manifolds? The difficulty lies in the inability to construct a reasonable test function on noncompact Finsler manifolds, and the more fundamental reason is that it is impossible to effectively control the second derivative of the distance function. This is due to the anisotropy at each point on the manifolds, hence the absence of a Lapalacian comparison theorem for the distance function on Finsler manifolds, even the manifolds are forward complete. As a nonlinear generalization of the weighted Riemannian space, Finsler metric measure space contains various non-Riemannian tensors. 
Such tensors are induced from the tangent bundle, forming the major obstacle to analysis on manifolds. The nonlinearity of the Finslerian Laplacian is totally different from the Riemannian one. As the references in this article named, nonlinear elliptic (or parabolic) operators on Riemannian manifolds or weighted Riemannian spaces often mean to add nonlinear terms to the Laplacian operator, or to replace the Laplacian by the $p$-Laplacian. However, the nonlinear operator on Finsler's metric measure space means that the Laplacian operator itself is nonlinear. Specifically, on Finsler manifolds, since derivation depends on the choice of direction vectors, the Finslerian Laplacian is intrinsically defined on the tangent bundle (or the pulled back bundle). When we project it onto the base manifold, the Laplacian operator on manifold depends on a local tangent vector field. For most studies, one can choose this local vector field as the gradient field of the function sought. However, this method fails when dealing with the comparison theorem of distance functions.

To improve the research, we first need to establish an appropriate and applicable Laplacian comparison theorem on the space of nonlinear measures of measures.
On the weighted Riemannian manifolds,
applying the curvature dimension condition $CD(K,N)$, G. Wei obtained the Laplacian comparison theorem and other results \cite{WW2009JDG}. K. Kuwae and X. Li proved a new comparison theorem with $N <1$ \cite{KL2021BLMS}. 
As a basic tool in geometric analysis, the Laplacian comparison theorem have been utilized widely on Riemannian geometric analysis, especially on the topic related to the elliptic and parabolic problems. However, most of these methods are far from complete in Finsler geometry.
As a nonlinear generalization of the weighted Riemannian space, a Finsler metric measure manifold admits a higher nonlinear Laplacian, so that the comparison theorems on Finsler manifolds are different and more difficult than the corresponding results on weighted Riemannian manifolds. More critically, application scope of a Finslerian comparison theorem is narrower, depending on the operator's linearization. Particularly, in Finsler cases, researchers provides different kinds of comparison theorems according to their requirements.
Hessian, Laplacian, and the volume comparison theorem on Finsler metric measure spaces was first generalized by B. Wu and Y. Xin in 2007 \cite{WuXin2007}. Later S. Ohta linked the curvature condition in the comparison theorem to the curvature-dimension condition in the sense of Lott-Sturm-Villani \cite{LV2009,St2006-1,St2006-2}. Based on it, he defined the weighted Ricci curvature and obtained the related comparison theorem\cite{Ohta2009}. Furthermore, Ohta and his collaborators have many applications in the study of Weitezenb\"ock's formula and heat flow. Meanwhile, B. Wu \cite{Wu2011PMD} and W. Zhao, Y. Shen \cite{ZS2013CJM} also studied different comparison theorems, respectively. As applications, They used the volume comparison theorem to obtain estimates of the fundamental group and the first Betti number, etc. On the other hand, S. Yin, Q. He and D. Zheng used the comparison theorem to obtain more properties of inequalities \cite{YHZ2014JIA}.\\

In this manuscript, We present a new Laplacian comparison theorem , which has more application prospects, by introducing mixed weighted Ricci curvature (cf. Definition \ref{Def-mwrc}), some non-Riemannian quantities (cf. \eqref{Def-U2}, \eqref{def-TVW}) and the uniform constant misalignment (cf. Definition \ref{def-misalignment} and Definition \ref{Def-Gmisalignment}). More precisely, we prove that
\begin{thm}\label{thm-LapComp-1}
	Let $(M,F,\mu,\alpha)$ be a forward complete $n$-dimensional Finsler metric measure space with finite misalignment $\alpha$.
	Denote the forward distance function by $r$ and by $V$ a fixed vector field on $M$. Suppose the mixed weighted Ricci curvature $^mRic^N_{\nabla r}$ of $M$ is bounded from below by $-K$ with $K>0$, for some $N>n$. Then, by setting $l=K/C(N,\alpha)$ with $C(N,\alpha)=N+(\alpha-1)n-\alpha$, wherever $r$ is smooth, the nonlinear Laplacian of $r$ with reference vector $V$ satisfies 
	\begin{eqnarray}
		\Delta^Vr\leq C(N,\alpha)\mathfrak{ct}_{-l}(r)+C_0,
	\end{eqnarray}
	where $C_0$ is a constant depending on $\alpha$ 
	and the bound of non-Riemannian curvatures $K_0$, namely, 
	$$F(U)+F^*(\mathcal{T})+F(\mathrm{div} C(V))\leq K_0,$$
	for any local vector field $V$. The definitions of $U$, $\mathcal{T}$ are given in \eqref{def-UyW} and \eqref{def-TVW}, respectively, and $\mathrm{div} C(V)=C^{ij}_{\,\,\,k\vert i}(V)V^k\frac{\delta}{\delta x^j}$.
\end{thm}

This comparison theorem also answer partially the question in \cite{BKT2024JFA} (Sec.6 III).

The boundedness of the misalignment $\alpha$ is equal to the boundedness of the uniform smoothness constant $\kappa$ and the uniform convexity constant $\kappa^*$ (cf. Theorem \ref{thm-kappa=alpha} in Section \ref{Ch-2}). The mixed weighted Ricci curvature and the non-Riemannian quantities $U$, $\mathcal{T}$ are first proposed in this article. More new concepts and basic properties of Finsler manifolds adopted in this theorem will be introduced in Section \ref{Ch-2}.

The conditions in this comparison theorem appear to be stronger than those used in the known comparison theorem. Especially, it contains more non-Riemannian quantities. But the most important feature of Theorem \ref{thm-LapComp-1} is that it can consistently overcome the difficulty of linearizing the Laplacian on other functions or tensors. \\

A direct example of the new Laplacian comparison theorem \ref{thm-LapComp-1} is the gradient estimation of the positive solution to different types of quasilinearized elliptic or parabolic equations on a noncompact forward complete Finsler manifold. 
In this manuscript, we consider the Finslerian Schr\"odinger equation, since it includes many other types of equations, such as the Poisson equation and the heat equation. On the other hand, it was Li-Yau's work on gradient estimation of the positive solutions to Schr\"odinger equation that named this geometric analysis method on Riemannian manifolds.

More precisely, the Finslerian Schr\"odinger equation studied in this paper can be seen as a generalization of this equation on more general metric measure spaces.
That is
\begin{eqnarray}\label{equ-FS1}
	(\Delta^{\nabla u}-\partial_t-q)u=0
\end{eqnarray}
which defined on a Finsler metric measure space $(M,F,\mu)$. Equation \eqref{equ-FS1} is more than a simple quasilinear parabolic equation. Its nonlinearity is also reflected in the fact that the space it located in is an anisotropic(asymmetric) metric space, geometrically, the metric depends on the tangent vector at each point. Such geometric nonlinarity fo the equation is named by ``Finslerian".

There can be many gradient estimates for positive solutions to Finslerian Schr\"odinger equation. We choose two typical gradient estimates to illustrate the idea, corresponding to two different locally weighted Ricci curvature conditions. Different weighted curvatures or mixed weighted curvatures can lead to different gradient estimates, which is also distinct from the Riemann case.

Unlike the cases on Riemannian manifolds, gradient estimates of Finslerian equations needs different curvature conditions on compact manifolds and noncompact manifolds, using the traditional methods. On compact manifolds, since there is no need to apply the comparison theorem, we can obtain the global gradient estimates by only requiring a lower bound of the weighted Ricci curvature. Specifically, we obtained Theorem \ref{thm-cpt-1} and Theorem \ref{thm-cpt-2} according to $Ric^N$ and $Ric^{\infty}$ bounded from below, respectively. Details are given in Section \ref{Ch-GGE}.

What we really care about is the situation of noncompact forward complete Finsler metric measure spaces. Now we need to require the lower boundedness of the mixed weighted Ricci curvature, as well as the bounded of misalignment and other non-Riemannian quantities, to ensure that the comparison theorem can be applied to control the second derivative of the distance function. Depending on the weighted Ricci curvature conditions, we obtain different gradient estimates. In the case with lower bounded mixed weighted Ricci curvature in direction $\nabla r$ with index $N>n$, i.e., $^mRic^N_{\nabla r}>-K$, we get the following theorem. 
\begin{thm}\label{thm-GEcomplete-1}
	Let $(M,F,\mu)$ be a forward complete Finsler metric measure space with or without boundary. Denote by $B(p,2R)$ the forward geodesic ball centered at $p$ with forward radius $2R$, which does not intersect the boundary $\partial M$ if it exists. Suppose the mixed weighted Ricci curvature $^mRic^N_{\nabla r}$ of $(M,F,\mu)$ admits a lower bound $-K(2R)$ in $B(p,2R)$ for a nonnegative constant $K(2R)$, the misalignment $\alpha$ is locally bounded from above by $A(2R)$, and the non-Riemannian tensors satisfy $F(U)+F^*(\mathcal{T})+F(\mathrm{div} C(V))\leq K_0$.
	Suppose equation \eqref{equ-FS1} satisfies that $q(x,t)$ is a function on $M\times [0,T]$ with $C^2$ and $C^1$ continuities in $x$ and $t$ variables, respectively. Moreover it satisfies $F(\nabla q(x))\leq\gamma(2R)$ and $\Delta^V q(x)\leq \theta(2R)$ in the distributional sense, for any point $x$ on $M$ and any local vector field $V$ around $x$, with two constants $\theta$ and $\gamma$ on $B(p,2R)$.
	
	Assume $u(x,t)$ is a positive solution on $M\times(0,T]$ to the Finslerian Schr\"odinger equation \eqref{equ-FS1}. It satisfies that
	\begin{eqnarray}\label{equ-thm1}
		\begin{split}
			\sup_{M\times [0,T]}&\left\{\frac{F^2(\nabla u)}{u^2}-\frac{\beta u_t}{u}-\beta q,\quad\frac{F^2(-\nabla u)}{u^2}-\frac{\beta u_t}{u}-\beta q\right\}\\
			\leq& C_3\beta^2R^{-2}(1+R+R\sqrt{K}+\frac{\beta^2}{\beta-1})+\frac{N}{2}\beta^2t^{-1}\\
			&+\left[\frac34\left(\frac{(\beta-1)^2\beta^8\gamma^4}{4\epsilon}\right)^{\frac13}+\frac{N^2\beta^4K^2}{4(1-\epsilon)(\beta-1)^2}+\frac{N}{2}\beta^3\theta\right]^{\frac12},
		\end{split}
	\end{eqnarray}
	for any $\beta>1$ and $\epsilon\in(0,1)$, where $C_3=C_3(N,A,K_0)$ is a constant depending on $N$, $A$ and $K_0$.
\end{thm}

Meanwhile, we get the following theorem when $^mRic^{\infty}_{\nabla r}$, the mixed weighted Ricci curvature with index $\infty$, is bounded from below locally. 
\begin{thm}\label{thm-GEcomplete-2}
	Let $(M,F,\mu)$ be a forward complete Finsler metric measure space with or without boundary. Denote by $B(p,2R)$ the forward geodesic ball centered at $p$ with forward radius $2R$, which does not intersect the boundary $\partial M$ if it exists. Suppose that the mixed weighted Ricci curvature $^mRic^{\infty}_{\nabla r}$ of $(M,F,\mu)$ admits a lower bound $-K(2R)$ and the $S$-curvature admits an $L^{\infty}$ upper bound $K'(2R)$ in $B(p,2R)$ for nonnegative constants $K(2R)$ and $K'(2R)$, the misalignment $\alpha$ is locally bounded from above by $A(2R)$, and the non-Riemannian tensors satisfy $F(U)+F^*(\mathcal{T})+F(\mathrm{div} C(V))\leq K_0$.
	Suppose equation \eqref{equ-FS1} satisfies that $q(x,t)$ is a function on $Q_{2R,T}:=B(p,2R)\times [t_0-T,t_0]$ with $C^2$ and $C^1$ continuities in $x$ and $t$ variables, respectively. 
	
	Assume $u(x,t)$ is a positive solution on $Q_{2R,T}$ to \eqref{equ-FS1}. Let $L$ be the upper bound of $u$ in $Q_{2R,T}$. It satisfies in the distributional sense that
	\begin{eqnarray}
		\begin{split}
			\sup_{Q_{2R,T}}\left\{\frac{F(\nabla u)}{u},\,\,\frac{F(-\nabla u)}{u}\right\}\leq &C(1+\log\frac{L}{u})\cdot\left[\frac{1}{R}+\frac{1}{\sqrt R}+\frac{1}{\sqrt T}+\sqrt K\right.\\
			&\quad\left.+\sqrt{K'}+\|q^-\|_{L^{\infty}(Q_{2R,T})}^{\frac12}+\|F(\nabla q)\|_{L^{\infty}(Q_{2R,T})}^{\frac13}\right],
		\end{split}
	\end{eqnarray}
	for some constant $C=C(N,A,K_0)$.
\end{thm}	

The proofs of Theorems \ref{thm-GEcomplete-1} and \ref{thm-GEcomplete-2} are given in Section \ref{Ch-noncpt}. 
Formally, our theorems on noncompact Finsler metric measure spaces require more and stricter conditions than the  theorems of gradient estimates in \cite{XiaC2014CVPDE} and \cite{Xia2023}, but if the Finsler metric is induced from a Riemannian oneb, both conditions in \cite{XiaC2014CVPDE}, \cite{Xia2023} and ours degenerate to the standard lower boundedness of the Ricci curvature or the weighted Ricci curvature. Methodologically, we need fewer analytical techniques, since we have more geometric requirements. Therefore, the approach proposed in this manuscript is a universal method for generalizing gradient estimation on noncompact forward complete Finsler metric measure spaces.\\

This manuscript is arranged as follows. In Section \ref{Ch-2}, we introduce some basic concepts and related properties on Finsler manifolds and metric measure spaces, as well as necessary theorems obtained by predecessors. In addition, we introduce misalignment and some new non-Riemann quantities in Subsection \ref{subsec-2.1}, and we introduce the mixed weighted Ricci curvature in Subsection \ref{subsec-2.2}. In Section \ref{sec-Lapcomp}, we give the proof of the new Laplacian comparison theorem, i.e., Theorem \ref{thm-LapComp-1}. We turn to introduce the fundamental properties of Finslerian Schr\"odinger's equation in Section \ref{SecFSE}, and then give gradient estimates on compact Finsler metric measure spaces in the traditional $CD(-K,N)$ conditions in Section \ref{Ch-GGE}. Precisely, we adopt the lower boundedness of the weighted Ricci curvature $Ric^{N}$ in Subsection \ref{Ch-GGE-1} to get a Li-Yau type gradient estimates with the $u_t$ term, and in Subsection \ref{Ch-GGE-2}, we provide Hamilton type gradient estimate without the $u_t$ term under the curvature condition about the weighted Ricci curvature $Ric^{\infty}$. In Section \ref{Ch-noncpt}, we consider gradient estimation of positive solutions on noncompact forward complete Finsler metric measure spaces and prove Theorem \ref{thm-GEcomplete-1} in curvature condition about $^mRic^N_{\nabla r}$, as well as some of its related corollaries. We also give the corresponding proof of Theorem \ref{thm-GEcomplete-2} with conditions about the mixed weighted curvature $^mRic^{\infty}_{\nabla r}$ in this section. Lastly, in Section \ref{Ch-ap}, we present applications of such gradient estimates, including Harnack inequality and several Liouville type theorems with the two different curvature conditions.

\section{Finsler metrics, tensors and Finsler metric measure spaces}\label{Ch-2}

In this section, we introduce basic concepts of Finsler metric measure spaces. According to the history, Finsler manifolds was proposed at the same period of Riemann geometry. This inspired the study of Finsler geometry from the perspective of metric geometry \cite{ShenShen2016}. To study the global geometric and topological properties of Finsler manifolds, researchers discovered several important non-Riemannian constants. Motivated by the generalization of the topological sphere theorem in Finsler geometry, H. B. Rademacher introduced the concepts of reversibility and reversible manifolds  \cite{Red2004}, which are also closely related to the analytical assumptions on Finsler manifolds. Later, K. Ball, E. Carlen and E. Lieb introduced the uniform smoothness and the uniform convexity in Banach space theory \cite{BCL1994Invent}, whose geometric explanation in Finsler geometry was given by S. Ohta \cite{Ohta2017}. We introduce a general concept called the misalignment, which contains the reversibility, uniform convexity and uniform smoothness constants. So we give the metric, connection and curvature properties in the first subsection. 

In the last two decades, encouraged by the developments of geometric measure theory and geometric analysis, Finsler geometry began to be revisited from a measure perspective. Such attempts could be traced back to the formulation of the $S$ curvature by Shen \cite{Shen1997}, to generalize the volume comparison theorem. Later, Ohta proposed the concept of Finsler metric measure space from the perspective of curvature-dimension condition $CD(K,N)$ in the sense of Lott-Sturm-Villani operators \cite{Ohta2009}\cite{Ohta2014}. We introduce a general concept called the mixed weighted Ricci curvature, which contains the weighted Ricci curvature as a special case. So we move on to the brief introduction of Finsler metric measure spaces in the second subsection.

\subsection{Basic tensors and important constants of Finsler manifolds}\label{subsec-2.1}

Let $M^n$ be a smooth, connected $n$-dimensional differential manifold. A \emph{Finsler metric} on $M$ is a generalization of Riemannian metric in the sense of norm on the tangent bundle. Briefly speaking, the \emph{Finsler structure} is a nonnegative function $F:TM\rightarrow [0,+\infty)$ satisfies that
\begin{itemize}
	\item[(i)] $F$ is smooth and positive on $TM\setminus\{0\}$;
	\item[(ii)] $F$ is a positive homogeneous norm, i.e., $F(x,ky)=kF(x,y)$ for any $(x,y)\in TM$ and for any $k>0$;
	\item[(iii)] $F$ is strongly pseudo-convex, namely, for any $(x,y)\in TM\setminus\{0\}$, the fundamental form matrix
	\begin{eqnarray}\label{Fiii}
		g_{ij}(x,y):=\frac12\frac{\partial F^2}{\partial y^i\partial y^j}(x,y)
	\end{eqnarray}
	is positive definite.
\end{itemize}

A Finsler manifold is a differential manifold endowed with a Finsler metric, which is always written as a pair $(M^n,F)$. The following basic knowledge about Finsler manifolds can be found in \cite{ShenShen2016}. The \emph{Cartan tensor} is a basic non-Riemannian tensor defined by
$$C(X,Y,Z):=C_{ijk}X^iY^jZ^k=\frac{1}{4}\frac{\partial^3F^2(x,y)}{\partial y^i\partial y^j\partial y^k}X^iY^jZ^k.$$ 
The \emph{mean cartan tensor} is $I=I_idx^i$, where locally $I_k=g^{jk}C_{ijk}$, which vanishes if and only if $F$ is deduced from a Riemannian metric due to the Deicke theorem.

There is a unique almost $g$-compatible and torsion-free connection on the pull back tangent bundle $\pi^*TM$ of the Finsler manifold $(M,F)$ called the \emph{Chern connection}. It is determined by
\begin{align*}
	\nabla_XY-\nabla_YX&=[X,Y];\\
	Z(g_y(X,Y))-g_y(\nabla_ZX,Y)-&g_y(X,\nabla_ZY)=2C_y(\nabla_Zy,X,Y),
\end{align*}
for any $X,Y,Z\in TM\setminus\{0\}$, where $C_y$ is the Cartan tensor. The Chern connection coefficients are locally denoted by $\Gamma^i_{jk}(x,y)$ in natural coordinate system, which could induce the \emph{spray coefficients} as $G^i=\frac12\Gamma^i_{jk}y^jy^k$.

The Chern connection can define the \emph{Chern Riemannian curvature} $R$ and \emph{Chern non-Riemannian connection} $P$. Denote by $\Omega$ the curvature form of Chern connection, so that
$$\Omega(X,Y)Z=R(X,Y)Z+P(X,\nabla_Yy,Z),$$
for any $X,Y,Z\in TM\setminus\{0\}$, where locally
$$R_{j\,\,kl}^{\,\,i}=\frac{\delta \Gamma^i_{jl}}{\delta x^k}+\frac{\delta \Gamma^i_{jk}}{\delta x^l}+\Gamma^i_{km}\Gamma^{m}_{jl}-\Gamma^i_{lm}\Gamma^{m}_{jk},\quad P_{j\,\,kl}^{\,\,i}=-\frac{\partial \Gamma^i_{jk}}{\partial y^l},$$
in which $\frac{\delta}{\delta x^i}=\frac{\partial}{\partial x^i}-N^j_i\frac{\partial}{\partial y^j}$, with $N^i_j=\frac{\partial G^i}{\partial y^j}$ called \emph{nonlinear connection}. 
We denote the horizontal Chern derivative by $``\vert"$ and the vertical Chern derivative by $``;"$. For example, 
$$v_{i\vert j}=\frac{\delta}{\delta x^j}v_i-\Gamma^k_{ij}v_k,\quad v_{i;j}=\frac{\partial }{\partial y^j}v_i,$$
for any 1- form $v=v_idx^i$ on the pull-back bundle.

The \emph{angular metric form} $h_y$ is defined by
\begin{eqnarray}\label{Def-amf}
	h_y(u, v)=g_y(u,v)-\frac{1}{F^2(y)}g_y(y, u)g_y(y, v),
\end{eqnarray}
for any $y, u, v\in T_xM$ with $y\neq 0$. Thus,
for any two linearly independent vector $y,v\in T_xM\setminus\{0\}$, which span a tangent plane $\Pi_y=\mathrm{span}\{y,v\}$, the \emph{flag curvature} with pole $y$ is defined by
\begin{eqnarray}\label{def-K}
	K(P,y) = K(y,u):=\frac{R_y(y, u, u, y)}{F^2(y)h_y(u, u)},
\end{eqnarray}
which is locally expressed by
$$K(y,u)=\frac{-R_{ijkl}(y)y^iu^jy^ku^l}{(g_{ik}(y)g_{jl}(y)-g_{il}(y)g_{jk}(y))y^iu^jy^ku^l}.$$

The \emph{Ricci curvature} is defined by
\begin{eqnarray}\label{def-Ric}
	Ric(y):=F^2(y)\sum_{i=1}^{n-1}K(y,e_i),
\end{eqnarray}
where $e_1,\cdots,e_{n-1},\frac{y}{F(y)}$ form an orthonormal basis of $T_xM$ with respect to $g_y$. On the other hand, Finsler non-Riemannian curvature tensors are abundant enough to make the geometry of (pull-back) tangent bundle colourful. The \emph{Landsberg curvature} of $(M,F)$ is given by
$$L:=L^i_{jk}\partial_i\otimes dx^j\otimes dx^k, \quad L^i_{jk}=-y^jP^{\,\,i}_{j\,\,kl}.$$
By the zero homogeniety, according to the Euler lemma, it is easy to see that 
$$C_{ijk}y^i=L_{ijk}y^i=0.$$

Noticing that the difference of two connections is a tensor. The \emph{$T$ curvature} (or \emph{tangent curvature}) introduced in Section 10.1 of \cite{ShenLec} is defined by
\begin{eqnarray}
	T_y(v):=g_y(D_vV,y)-\hat g_x(\hat D_vV,y),
\end{eqnarray}
where $v\in T_xM$, $V$ is a vector field with $V(x)=v$, and $\hat D$ denotes the Levi-Civita connection of the induced Riemannian metric $\hat g=g_Y$. The $T$ curvature vanishes if and only if the Chern non-Riemannian curvature $P$ vanishes. Similarly, let $\{e_i\}$ be an orthonormal basis with respect to metric $g(x,V)$ at point $x$, where $V$ is a fixed reference vector field. Moreover, Let $\{E_i\}$ be the local vector fields obtained by moving $\{e_i\}$ in parallel in a neighborhood of $x$ on $M$. We define tensor $U$ by
\begin{eqnarray}\label{Def-U3}
	U(V,W,Z)=g(x,Z)\left(\sum_{i=1}^n (D^W_{e_i}E_i-D^V_{e_i}E_i),Z\right),
\end{eqnarray}
where $V,W,Z$ are three vector fields on $M$, and $\{e_i\}$, $\{E_i\}$ form an orthonormal frame and induced locally parallel vector fields, respectively, with respect to the metric $g(x,V)$. In paticular, we can define a tensor $U$ from \eqref{Def-U3} by
\begin{eqnarray}\label{Def-U2}
	U_y(W)=U(Y,W,W)=g(x,W)\left(\sum_{i=1}^n (D^W_{e_i}E_i-\hat D_{e_i}E_i),W\right),
\end{eqnarray}
for any local vector field $W$, where $Y$ is the geodesic extension field of $y$ in a neighborhood of $x$, and $\hat D$ denotes the Levi-Civita connection of the induced Riemannian metric $\hat g=g_Y$. 
We also denote $U_y(W)=g(x,W)(U(y,W),W)$, with 
\begin{eqnarray}\label{def-UyW}
	U(y,W)=\sum_{i=1}^n (D^W_{e_i}E_i-\hat D_{e_i}E_i)
\end{eqnarray}
a vector field on the sphere bundle.
It is also a curvature, which may be considered as a kind of trace of the $T$ curvature.

For any given non-vanishing vector field $V$, $g_{ij}(x,V)$ could induce a Riemannian structure $g_V$ on $M$ via
\begin{eqnarray*}
	g_V(X,Y)=g_{ij}(x,V)X^iY^j,
\end{eqnarray*}
for any vector pairs $X,Y\in T_xM$.
The norm of $V$ is measured as $F^2(V)=g_V(V,V)$. 

The Riemannian structure is inconsistent when the reference vectors are different. For example, given three different local non-vanishing vector fields around $x$, namely, $V,W,Y$, the norm of $Y$ about $g_V$ and $g_W$ maybe not the same in general case. The ratio $g_V(Y,Y)/g_W(Y,Y)$ is a function about $V,W,Y$. We give the following definition.
\begin{defn}[Local misalignment]\label{def-misalignment}
	For a Finsler manifold $(M,F)$, we define the local misalignments of the Finsler metric with respect to the local nonvanishing vector field $W$ by 
	\begin{align}
		\alpha_M(x,W)&=\sup_{V\in S_xM}\sup_{Y\neq0}\frac{g_V(Y,Y)}{g_W(Y,Y)}=\sup_{V\in S_xM;\,g_W(Y,Y)=1}g_V(Y,Y),\\
		\alpha_m(x,W)&=\inf_{V\in S_xM}\inf_{Y\neq0}\frac{g_V(Y,Y)}{g_W(Y,Y)}=\inf_{V\in S_xM;\,g_W(Y,Y)=1}g_V(Y,Y),
	\end{align}
	where $\alpha_M(x,W)$ is called the local maximal misalignment and $\alpha_m(x,W)$ is called the minimal misalignment at $x$ with respect to $W$. Moreover, we define the local misalignments by
	\begin{align}
		\alpha_M(x)&=\sup_{W\neq 0}\alpha_M(x,W)=\sup_{V,W,Y\in S_xM}\frac{g_V(Y,Y)}{g_W(Y,Y)},\\
		\alpha_m(x)&=\inf_{W\neq 0}\alpha_m(x,W)=\inf_{V,W,Y\in S_xM}\frac{g_V(Y,Y)}{g_W(Y,Y)}.
	\end{align} 
\end{defn}

\begin{rem}
	According to the roles of the vector fields $V$ and $W$ in the definition above, the local misalignments of the Finsler metric satisfies that 
	\begin{eqnarray}
		\alpha_M(x)\alpha_m(x)=1,
	\end{eqnarray}
	for any $x$. So when we mention the local misalignment of a Finsler metric, we always refer to $\alpha_M(x)$. That is, we denote $\alpha(x)=\alpha_M(x)$ locally.
\end{rem}

\begin{defn}[Global misalignment]\label{Def-Gmisalignment}
	For a Finsler manifold $(M,F)$, the global misalignments of a Finsler metric with respect to a global nonvanishing vector field $W$ are defined by
	\begin{eqnarray}
		\alpha_M(W)=\sup_{x\in M}\alpha_M(x,W)\quad\mbox{and}\quad \alpha_m(W)=\inf_{x\in M}\alpha_m(x,W).
	\end{eqnarray}
	Moreover, the global misalignment of the Finsler metric is defined by
	\begin{eqnarray}\label{def-alpha}
		\alpha=\sup_{x\in M}\alpha_M(x)=\sup_{x\in M}\sup_{V,W,Y\in S_xM}\frac{g_V(Y,Y)}{g_W(Y,Y)}.
	\end{eqnarray}
\end{defn}

\begin{rem}
	\begin{itemize}
		\item[(i)] One would also easily find that 
		$$\frac{1}{\alpha}\leq\alpha_m(W)\leq1\leq\alpha_M(W)\leq\alpha,$$
		for any local vector field $W$. So given any vector field $W$, there exists a positive constant $C\geq1$ such that $0< \frac1C\leq\alpha_m(W)\leq 1\leq\alpha_M(W)\leq C$, provided that $\alpha\leq C$.
		\item[(ii)] The geometric meaning of the local misalignment is the maximal deviation of the unit sphere with respect to the metric $g_W$ and any other $g_V$ with another reference vector field $V$ at point $x$. So we also can use 
		\begin{eqnarray*}
			\alpha(x,W)=\max\{\alpha_M(x,W),\frac{1}{\alpha_m(x,W)}\}\quad\mbox{and}\quad\alpha(W)=\max\{\alpha_M(W),\frac{1}{\alpha_m(W)}\}
		\end{eqnarray*}
		to discribe the misalignment with respect to $W$.
	\end{itemize}
\end{rem}

When the Finsler metric is induced from a Riemannian one, the misalignment is obviously equal to 1, since $g(x,y)=g(x)$ in this case. On the other hand, if a Finsler metric $F$ satisfies that $\alpha_M=1$, it means at any point $x$, the norm of a vector $Y$ is coincident with respect to any induced Riemannian structure $g_V$, for any $V$. That is, the Finsler fundamental form $g(x,y)$ is independent of $y$. Consequently, we arrive at
\begin{prop}
	A Finsler manifold $(M,F)$ is a Riemannian manifold if and only if $\alpha_M=1$.
\end{prop}
Since that, we give an important class of Finsler manifold as the following.
\begin{defn}
	We call a Finsler manifold $(M,F)$ has finite misalignment if $\alpha<+\infty$.
\end{defn}



A Finsler metric is defined to be \emph{reversible} if $F(x,V)=\bar F(x,V)$ for any point $x$ and any vector field $V$, where $\bar F(x,V):=F(x,-V)$ is called the \emph{reversed Finsler metric} of $F$. We define the \emph{reversibility} of $(M,F)$ by
\begin{eqnarray*}
	\rho:=\sup_{x\in M}\sup_{V\neq 0}\frac{F(x,V)}{\bar F(x,V)}.
\end{eqnarray*}
Obviously, $F$ is \emph{reversible} if and only if $\rho\equiv 1$. A Finsler manifold $(M,F)$ is said to have \emph{finite reversibility} if $\rho<+\infty$.

It seems that the reversibility describes the difference of norms with respect to the Finsler metric along geodesics and the misalignment describes those difference on the geodesic sphere. However, we will soon know that the misalignment is a more general concept, that is, the reversibility is a special case of it.

Moreover, we say $F$ satisfies \emph{uniform convexity} and \emph{uniform smoothness} if there exist uniform positive constants $\kappa^*$ and $\kappa$, called the \emph{uniform convexity constant} and \emph{uniform smoothness constant}, respectively, such that for any $x\in M$, $V\in T_xM\setminus\{0\}$ and $y\in T_xM$, we have 
\begin{eqnarray}\label{kappa}
	\kappa^*F^2(x,y)\leq g_V(y,y)\leq\kappa F^2(x,y),
\end{eqnarray}
where $g_V=(g_{ij}(x,V))$ is the Riemannian metric on $M$ induced from $F$ with respect to the reference vector $V$. In this situation, the reversibility $\rho$ could be controlled by $\kappa$ and $\kappa^*$ as
\begin{eqnarray}
	1\leq \rho\leq \min\{\sqrt{\kappa},\sqrt{1/\kappa^*}\}.
\end{eqnarray}
$F$ is Riemannian if and only if $\kappa=1$ if and only if $\kappa^*=1$ \cite{Ohta2017}.

The bounds of misalignment $\alpha$ provides the bounds of the other constants, including $\rho$, $\kappa^*$ and $\kappa$. In fact, we have 
\begin{eqnarray}
	\rho^2=\sup_{x\in M}\sup_{V\in S_xM}\frac{g_V(V,V)}{g_{-V}(V,V)}\leq \alpha,
\end{eqnarray}
by choosing $Y=V=V$ and $W=-V$ in \eqref{def-alpha}, the definition of $\alpha$. Furthermore, 
\begin{eqnarray}
	\kappa\geq \frac{g_V(Y,Y)}{g_Y(Y,Y)}\geq \frac{1}{\alpha},
\end{eqnarray}
and we have actually from the definition of $\kappa$ that
\begin{eqnarray}
	\kappa=\sup_{x\in M}\sup_{V,Y\in S_xM}\frac{g_V(Y,Y)}{g_Y(Y,Y)}\leq\alpha, 
\end{eqnarray}
by choosing $V=V$ and $W=Y=Y$ in \eqref{def-alpha}. Thus, it arrives at  $1/\alpha\leq \kappa\leq \alpha$. Meanwhile, the same process shows $\frac{1}{\alpha}\leq\kappa^*\leq \alpha$. Now we have proved that 
\begin{prop}\label{propkapparhoalpha}
	A Finsler manifold $(M,F)$ satisfies finite reversibility, uniform convexity and uniform smoothness if it has finite misalignment.
\end{prop}

On the other hand, we can see from the definition of the uniform convexity and the uniform smoothness that
\begin{thm}\label{thm-kappa=alpha}
	A Finsler manifold $(M,F)$ is uniform convexity and uniform smoothness if and only if it satisfies finite misalignment.
\end{thm}
\begin{proof}
	Suppose $(M,F)$ satisfies \eqref{kappa} for any $x\in M$, $V\in T_xM\setminus\{0\}$ and $y\in T_xM$. Then it follows from \eqref{def-alpha} that
	\begin{eqnarray}
		\alpha=\sup_{x\in M}\sup_{V,W,Y\in S_xM}\frac{g_V(Y,Y)}{g_W(Y,Y)}\leq\sup_{x\in M}\sup_{Y\in S_xM}\frac{\kappa F^2(x,Y)}{\kappa^*F^2(x,Y)}=\frac{\kappa}{\kappa^*}.
	\end{eqnarray}
	The sufficiency is given in Proposition \ref{propkapparhoalpha}.
\end{proof}

\noindent

At last, we give some examples of Finsler manifolds with finite or infinite misalignment.
\begin{exm}\label{example}
	\begin{itemize}
		\item[i)] Riemannian manifolds are trivial examples of Finsler manifolds with finite misalignment.
		\item[ii)] A Minkowski manifold is a Finsler manifold with finite misalignment.
		\item[iii)] A unit disk with the Funk metric is a Finsler manifold with infinity misalignment.
	\end{itemize} 
\end{exm}


\begin{rem}\label{rem-localalpha}
	\begin{itemize}
		\item[i)] According to Theorem \ref{thm-kappa=alpha}, the condition of finite misalignment through this article can be replaced alternatively by the uniform convexity and uniform smoothness conditions.
		\item[ii)] The third example shows the existence of noncompact forward complete Finsler manifolds with infinite misalignment, which limits the scope of our application of this approach. However, referring to the ideal in \cite{KZ2022JMPA}, we can assume that in the noncompact case, the manifold has a locally finite misalignment. That is, on any forward bounded closed set $\Omega$ of the Finsler manifold, there is a constant $A$ dependent on $\Omega$ such that $\alpha(x)$ is less than $A$ everywhere on $\Omega$.
	\end{itemize}	
\end{rem}

\subsection{Nonlinear Laplacian and weighted Ricci curvature on Finsler metric measure spaces}\label{subsec-2.2}

Noticing that $F$ is actually a norm on the tangent bundle $TM$, there is a natural \emph{dual norm} $F^*$ on the cotangent bundle $T^*M$ defined by
\begin{eqnarray*}
	F^*(x,\xi):=\sup_{F(x,y)=1}\xi(y),
\end{eqnarray*}
for any $\xi \in T_x^*M$. The fundamental form of $F^*$ is defined by
\begin{eqnarray*}
	g^{*ij}(x,\xi):=\frac12\frac{\partial^2F^{*2}}{\partial\xi_i\partial \xi_j},
\end{eqnarray*}
which is positive definite for each $(x,\xi)\in T^*M\setminus\{0\}$ and is the inverse matrix of $g_{ij}(x,y)$ defined by (\ref{Fiii}).

The \emph{Legendre transformation} is an isomorphism between $T_xM$ and $T_x^*M$ given by a map $l$ with
\begin{eqnarray*}
	l(y):=\begin{cases}
		g_y(y,\cdot)\,\,\quad \,\mbox{for } y\in T_xM\setminus\{0\},\\
		0 \quad\quad\quad\quad\mbox{for } y=0.
	\end{cases}
\end{eqnarray*}

One can verify that $g^*_{ij}(x,\xi)=g^{ij}(x,y)$, for $\xi=l(y)$. If $F$ is uniformly smooth or uniformly convex, then so is $F^*$. Precisely, there exist two constants $\tilde\kappa=(\kappa^*)^{-1}$ and $\tilde\kappa^*=\kappa^{-1}$ such that for any $x\in M$, $\xi\in T_x^*M\setminus\{0\}$ and $\zeta\in T^*_xM$, 
\begin{eqnarray}
	\tilde\kappa^*F^{*2}(x,\zeta)\leq g_{\xi}^*(\zeta,\zeta)=g^{ij}(x,l^{-1}(\xi))\zeta_i\zeta_j\leq\tilde\kappa F^{*2}(x,\zeta),
\end{eqnarray}
if \eqref{kappa} holds.

For any smooth function $f:M\rightarrow \mathbb{R}$, $df$ denotes its differential 1-form and its \emph{gradient} $\nabla f$ is defined as the dual of the 1-form via the Legendre transformation, namely, $\nabla f(x):=l^{-1}(df(x))\in T_xM$. Locally it can be written as
\begin{eqnarray*}
	\nabla f=g^{ij}(x,\nabla f)\frac{\partial f}{\partial x^i}\frac{\partial }{\partial x^j}
\end{eqnarray*}
on $M_f:=\{df\neq 0\}$. The Hessian of $f$ is defined via the Chern connection by
\begin{eqnarray*}
	\nabla^2f(X,Y)=g_{\nabla f}(\nabla_X^{\nabla f}\nabla f,Y).
\end{eqnarray*}
It can be shown that $\nabla^2f(X,Y)$ is symmetric \cite{Ohta2014}.\\

For any two points $p,q$ on $M$, the \emph{distance function} is defined by
$$d_p(q):=d(p,q):=\inf_{\gamma}\int_0^1F(\gamma(t),\dot\gamma(t))dt,$$
where the infimum is taken over all the $C^1$ curves $\gamma:[0,1]\rightarrow M$ such that $\gamma(0)=p$ and $\gamma(1)=q$. 
Fixed a base point $p$ on $M$, we denote the forward distance function by $r$. That is, $r(x)=d(p,x)$, with $d$ denotes the forward distance. 

The \emph{forward distance function} $r$ is a function defined on the Finsler manifold $M$. $dr$ is a 1-form on $M$, whose dual is a gradient vector field, noted by $\nabla r$. Precisely, $\nabla r=g^{ij}(x,\nabla r)\frac{\partial r}{\partial x^i}\frac{\partial}{\partial x^j}$. Taking the Chern horizontal derivative of $\nabla r$ yields the\emph{ Hessian of distance function} $H(r)=\nabla^2r$. Locally, in natural coordinate system 
\begin{eqnarray}\label{Hessianr}
	\nabla^2 r(\frac{\partial}{\partial x^i},\frac{\partial }{\partial x^j})=\nabla_j\nabla_ir=\frac{\partial^2 r}{\partial x^i\partial x^j}-\Gamma^k_{ij}(x,\nabla r)\frac{\partial r}{\partial x^k}.
\end{eqnarray} 
In \eqref{Hessianr}, the derivative is taken in the direction $\nabla r$ naturally. Generally, we can take the derivative in any direction. Suppose $V$ is a local vector field around $x$ on $M$. We define the mixed Hessian of a function $f$ on $M$ by $H^V(f)=\nabla^{V2} f$. Locally
\begin{eqnarray}\label{mixedHessianf}
	H_{ij}^V(f)=\nabla^{V2} f(\frac{\partial}{\partial x^i},\frac{\partial }{\partial x^j})=\nabla^V_j\nabla^V_if=\frac{\partial^2 f}{\partial x^i\partial x^j}-\Gamma^k_{ij}(x,V)\frac{\partial f}{\partial x^k}.
\end{eqnarray}



Note that the distance function may not be symmetric about $p$ and $q$ unless $F$ is reversible. A $C^2$ curve $\gamma$ is called a \emph{geodesic} if locally
$$\ddot\gamma(t)+2G^i(\gamma(t),\dot \gamma(t))=0,$$
where $G^i(x,y)$ are the \emph{spray coefficients}.
A \emph{forward geodesic ball} centered at $p$ with radius $R$ can be represented by
\begin{eqnarray*}
	B_R^+(p):=\{q\in M\,:\, d(p,q)<R\}.
\end{eqnarray*}
Adopting the exponential map, a Finsler manifold $(M,F)$ is said to be \emph{forward complete} or forward geodesically complete if the exponential map is defined on the entire $TM$. Thus, any two points in a forward complete manifold $M$ can be connected by a minimal forward geodesic and the forward closed balls $\overline{B_R^+(p)}$ are compact.\\

A \emph{Finsler metric measure space} $(M,F,d\mu)$ is a Finsler manifold equipped with an given measure $\mu$. 
In local coordinates $\{x^i\}_{i=1}^n$, we can express the volume form as $d\mu=\sigma(x)dx^1\wedge\cdots\wedge dx^n$ with $\sigma(x)>0$. For any $y\in T_xM\setminus\{0\}$, define 
$$\tau(x,y):=\log\frac{\sqrt{\det g_{ij}(x,y)}}{\sigma(x)},$$
which is called the \emph{distortion} of $(M,F,\mu)$. The definition of the \emph{S-curvature} is given in the following.
\begin{defn}\cite{Shen1997}\cite{ShenShen2016}\label{def-S}
	Suppose $(M,F,d\mu)$ is a Finsler metric measure space. For any point $x\in M$, let $\gamma=\gamma(t)$ be a forward geodesic from $x$ with the initial tangent vector $\dot\gamma(0)=y$. The S-curvature of $(M,F,d\mu)$ is
	\begin{eqnarray*}
		S(x,y):=\frac{d}{dt}\tau=\frac{d}{dt}(\frac12\log\det(g_{ij})-\log\sigma(x))(\gamma(t),\dot\gamma(t))\mid_{t=0}.
	\end{eqnarray*}
\end{defn}
Definition \ref{def-S} means that the S-curvature is the changing of distortion along the geodesic in direction $y$. We define the \emph{discrete difference of} $\nabla\tau$ on the tangent sphere, denoted by $\mathcal{T}$, as
\begin{eqnarray}\label{def-TVW}
	\mathcal{T}(V,W):=\nabla^V\tau(V)-\nabla^W\tau(W),
\end{eqnarray}
where $\nabla^V$ is the Chern connectionwith reference vector field $V$. Its local expression is $\mathcal{T}(V,W)=\mathcal{T}_i(V,W)dx^i$, given vector fields $V,W$ on $M$, with 
\begin{eqnarray}
	\mathcal{T}_i(V,W)=\frac{\delta}{\delta x^i}\tau(V)-\frac{\delta}{\delta x^i}\tau(W),
\end{eqnarray}
where $\frac{\delta}{\delta x^k}\tau=\frac{\partial}{\partial x^k}\tau-N^j_k\frac{\partial}{\partial y^j}\tau$ could also be denoted by $\tau_i$. Obviously, $\mathcal{T}$ is anti-symmetric about $V$ and $W$. 

If in local coordinates $\{x^i\}_{i=1}^n$, expressing $d\mu=e^{\Phi}dx^1\cdots dx^n$, the \emph{divergence} of a smooth vector field $V$ can be written as
\begin{eqnarray}
	\mathrm{div}_{d\mu}V=\sum_{i=1}^n(\frac{\partial V^i}{\partial x^i}+V^i\frac{\partial \Phi}{\partial x^i}).
\end{eqnarray}
The \emph{Finsler Laplacian} of a function $f$ on $M$ could now be given by
\begin{eqnarray}\label{def-WLap}
	\Delta_{d\mu}f:=\mathrm{div}_{d\mu}(\nabla f).
\end{eqnarray}
Noticing that $\Delta_{d\mu}f=\Delta_{d\mu}^{\nabla f}f$, where $\Delta_{d\mu}^{\nabla f}f:=\mathrm{div}_{d\mu}(\nabla^{\nabla f}f)$ is in the view of weighted Laplacian with
\begin{eqnarray}
	\nabla^{\nabla f}f:=\begin{cases}
		&g^{ij}(x,\nabla f)\frac{\partial f}{\partial x^i}\frac{\partial}{\partial x^j} \quad \mbox{for } x\in M_f;\\
		&0\quad\quad\quad\quad\quad\quad\quad\,\,\, \mbox{for }x\notin M_f.
	\end{cases}
\end{eqnarray}

On a Finsler metric measure space $(M,F,d\mu)$, we always denote $\Delta_{d\mu}^{\nabla f}=\Delta f$ and $\nabla ^{\nabla f}f=\nabla f$ for short unless otherwise indicated. On the other hand, the Laplacian of a function $f$ on a Riemannian manifold is the trace of the Hessian of $f$ with respect to the Riemannian metric $g$. On $(M,F,d\mu)$, the \emph{weighted Hessian} $\tilde{H}(f)$ of a function $f$ on $M_f=\{x\in M:\nabla f\mid_x \neq 0\}$ is defined in \cite{Wu2015} by
\begin{eqnarray}\label{Def-WHf}
	\tilde{H}(f)(X,Y)=H(f)(X,Y)-\frac{S(\nabla f)}{n-1}h_{\nabla f}(X,Y),
\end{eqnarray}
where $h_{\nabla f}$ is the angular metric form in the direction $\nabla f$, given in \eqref{Def-amf} and $H(f)(X,Y)=\nabla^2(f)(X,Y)=XY(f)-\nabla_X^{\nabla f} Y(f)$ is the usual Hessian of function $f$. It is clear that $\tilde{H}(f)$ is still a symmetric bilinear form with
\begin{eqnarray}
	\Delta f=\mathrm{tr}_{\nabla f}\tilde{H}(f).
\end{eqnarray}

Inspired by this, we define the \emph{mixed weighted Hessian} \begin{eqnarray}\label{Def-MWH}
	\tilde{H}^V(f)(X,Y)
	:=\nabla^2(f)(X,Y)-\frac{S(\nabla f)}{n-1}h_{V,\nabla f}(X,Y),
\end{eqnarray}
where $h_{V,\nabla f}$ is the \emph{mixed angular metric form in the directions $V$ and $\nabla f$}, which is defined by
\begin{eqnarray}
	h_{V,\nabla f}(X,Y)=g_V(X,Y)-\frac{1}{F_V^2(\nabla f)}g_V(X,\nabla f)g_V(Y,\nabla f),
\end{eqnarray}
for any vector $X,Y$.

It is necessary to remark that $h_{\nabla f,\nabla f}=h_{\nabla f}$, so that $\tilde H^{\nabla f}(f)=\tilde H(f)$ for any function $f$ on $M$.\\



A Finsler Laplacian is better to be viewed in a weak sense due to the lack of regularity. Concretely, assuming $f\in W^{1,p}(M)$,
\begin{eqnarray*}
	\int_M\phi\Delta_{d\mu}fd\mu=-\int_Md\phi(\nabla f)d\mu,
\end{eqnarray*}
for any test function $\phi\in C^{\infty}_0(M)$.

With the assistance of the S-curvature, one can present the definition of the\emph{ weighted Ricci curvature} as the following.	
\begin{defn}\cite{Ohta2014}\cite{ShenShen2016}
	Given a unit vector $V\in T_xM$ and an positive number $k$, the weighted Ricci curvature is defined by
	\begin{eqnarray*}
		Ric^k(V):=\begin{cases}
			&Ric(x,V)+\dot{S}(x,V)\quad\mbox{if } S(x,V)=0 \mbox{ and } k=n \mbox{ or if }k=\infty;\\
			&-\infty\quad\quad\quad\quad\quad\quad\quad\quad\quad\quad\quad\quad\mbox{if }S(x,V)\neq0 \mbox{ and if } k=n;\\
			&Ric(x,V)+\dot{S}(x,V)-\frac{S^2(x,V)}{k-n}\quad\quad\quad\quad\quad\quad\quad\mbox{ if }n<k<\infty,
		\end{cases}
	\end{eqnarray*}
	where the $Ric$ is defined in \eqref{def-Ric}, and the derivative ``$\dot{\quad}$" is taken along the geodesic started from $x$ in the direction of $V$.
\end{defn}
According to the definition of weighted Ricci curvature, B. Wu defined the \emph{weighted flag curvature} when $k=N\in (1,n)\cup(n,\infty)$ in \cite{Wu2015}. We introduce this concept for any $k$ in the following.
\begin{defn}[weighted flag curvature]\label{Def-wfc}
	Let $(M,F,d\mu)$ be a Finsler metric measure space,
	and $V,W\in T_xM$ be linearly independent vectors. The
	weighted flag curvature $K^k(V; W)$ is defined by
	\begin{eqnarray*}
		K^k(V; W):=\begin{cases}
			&K(V; W)+\frac{\dot S(V)}{(n-1)F^2(V)}\,\,\,\mbox{if } S(x,V)=0 \mbox{ and } k=n \mbox{ or if }k=\infty;\\
			&-\infty\quad\quad\quad\quad\quad\quad\quad\quad\quad\quad\quad\quad\,\,\mbox{if }S(x,V)\neq0 \mbox{ and if } k=n;\\
			& K(V; W)+\frac{\dot S(V)}{(n-1)F^2(V)}-\frac{S^2(V)}{(n-1)(k-n)F^2(V)}
			
			\quad\quad\quad\mbox{ if }n<k<\infty,
		\end{cases}
	\end{eqnarray*}
	where $K$ is defined in \eqref{def-K}, and the derivative ``$\dot{\quad}$"  is taken along the geodesic started from $x$ in the direction of $V$.
\end{defn} 
Moreover, we define the \emph{mixed weighted Ricci curvature}, denoted by $^mRic^{k}$.
\begin{defn}[mixed weighted Ricci curvature]\label{Def-mwrc}
	Given two unit vectors $V,W\in T_xM$ and an positive number $k$, the mixed weighted Ricci curvature $^mRic^k(V,W)=\,^mRic^k_W(V)$ is defined by
	\begin{eqnarray*}
		^mRic^k_{W}(V):=\begin{cases}
			&\mathrm{tr}_WR_{V}(V)+\dot{S}(x,V)\quad\mbox{if } S(x,V)=0 \mbox{ and } k=n \mbox{ or if }k=\infty;\\
			&-\infty\quad\quad\quad\quad\quad\quad\quad\quad\quad\quad\quad\quad\mbox{if }S(x,V)\neq0 \mbox{ and if } k=n;\\
			&\mathrm{tr}_WR_{V}(V)+\dot{S}(x,V)-\frac{S^2(x,V)}{k-n}\quad\quad\quad\quad\quad\quad\quad\mbox{ if }n<k<\infty,
		\end{cases}
	\end{eqnarray*}
	where the derivative is taken along the geodesic started from $x$ in the direction of $V$, and $\mathrm{tr}_WR_{V}(V)=g^{ij}(W)g_{V}(R_V(e_i,V)V,e_j)$ means taking trace of the flag curvature with respect to $g(x,W)$.
\end{defn}

\begin{rem}
It looks like the mixed weighted curvature $^m\operatorname{Ric}^k_W(V)$ is defined only at the point where the vector field $W$ is non-zero. However, we can expand it when $W$ does not identically vanish. For a locally nontrivial vector field $W$, with the set $M_W:=\{x\in M \mid W(x)\neq0\}$ has zero measure. We can choose a class of orthonormal frame fields $\{e_1,\cdots,e_{n-1},e_n=\frac{V}{F_W(V)}\}$ with respect to $g_W=g_{ij}(x,W)dx^i\otimes dx^j$ wherever $W\neq 0$, so that the Ricci curvature part in Definition \ref{Def-mwrc} is given by
\begin{eqnarray}\label{mRicw}
	\operatorname{tr}_WR_V(V):=F_W^2(V) \sum_{i=1}^{n-1} K\left(e_i,V \right),
\end{eqnarray}
where $F_W$ is the norm with reference vector $W$. For the point where $W=0$, we take a limit on the orthogonal frame, and define the limitation of (\ref{mRicw}) as the value of $\operatorname{tr}_WR_V(V)$ in the mixed weighted Ricci curvature at $W=0$. Because the definition of $\operatorname{tr}_WR_V(V)$ as a tensor is independent of the choice of the frames, so this limit is well defined. Moreover, Since the indicatrix at each point of $M$ is compact for a regular Finsler metric, we always have  
$$\frac{1}{\alpha(x)}F^2(V)\leq F_W^2(V)\leq \alpha(x)F^2(V).$$
Therefore, we define the mixed weighted Ricci curvature $^mRic^k_{W}(V)$ to be the weighted Ricci curvature $Ric^k(V)$ when the reference vector $W=0$.
\end{rem}
\begin{rem}
	\begin{itemize}
		\item[i)]	The weighted Ricci curvature is a special case of the mixed weighted Ricci curvature, i.e., $Ric^k(V)=\,^mRic^k_V(V)$.
		\item[ii)] We call the mixed weighted Ricci curvature bounded from below by a constant $K$, denoted by $^m\operatorname{Ric}^k_W\geq K$, if the inequality holds that $^m\operatorname{Ric}^k_W(V)\geq KF^2(V)$.
	\end{itemize} 
\end{rem}

In \cite{Wu2015}, Wu adopted the weighted flag curvature to obtain the weighted Hessian comparison theorem in the help of the function $\mathfrak{ct}_c(r)$, i.e., 
\begin{eqnarray}
	\mathfrak{ct}_c(r)=\begin{cases}
		\sqrt{c}\cot\sqrt{c}r,\quad\,\,\,\quad c>0,\\
		1/r, \quad\quad\quad\quad\quad\,\,\,\,  c=0,\\
		\sqrt{-c}\coth
		\sqrt{-c}r, \quad c<0.
	\end{cases}
\end{eqnarray}

\begin{thm}[\cite{Wu2015}]\label{thm-HessComp}
	Let $(M,F,d\mu)$ be a Finsler metric measure space and $r=d_F(p,\cdot)$ be the forward distance from $p$. Suppose that for
	some $N>n$, the weighted flag curvature of $M$ satisfies $K^N\geq\frac{N-1}{n-1}c$, where $c$ is a positive constant. Then, for
	any vector field $X$ on $M$, the following inequality is valid whenever $r$ is smooth. Particularly, $r < \pi/\sqrt{c}$ when $c>0$.
	\begin{eqnarray}
		\tilde H(r)(X,X)\leq \frac{N-1}{n-1} \mathfrak{ct}_c(r)(g_{\nabla r}(X,X)-g_{\nabla r}(\nabla r,X)^2).
	\end{eqnarray}
\end{thm}

\section{A new kind of Laplacian Comparison Theorem}\label{sec-Lapcomp}


The main content in this section is to establish the new comparison theorem \ref{thm-LapComp-1} on Finsler metric measure spaces, whose curvature conditions includes some purely non-Riemannian curvatures.

Suppose $(M,F,\mu)$ is a forward complete Finsler metric measure space. We denote by $V$ a fixed vector field on $M$. The nonlinear Laplacian of a function $f$ on $M$ with respect to $V$ is given by
\begin{eqnarray}\label{Lap-r-1}
	\Delta^Vf=e^{-\Phi}\frac{\partial}{\partial x^i}(e^{\Phi}g^{ij}(x,V)\frac{\partial f}{\partial x^j})=\frac{\partial \Phi}{\partial x^i}g^{ij}(x,V)\frac{\partial f}{\partial x^j}+\frac{\partial}{\partial x^i}(g^{ij}(x,V)\frac{\partial f}{\partial x^j}).
\end{eqnarray}
We always denote the partial derivative of a function by its lower indices, such as $\Phi_i=\frac{\partial \Phi}{\partial x^i}$, for short. We want to establish a Laplacian comparison theorem for such nonlinear elliptic operator. 

With the assistance of the misalignment, one can get the following Theorem \ref{thm-LapComp-2} from \eqref{Lap-r-4} according to
Theorem \ref{thm-HessComp} and the curvature conditions. It provides a version of the comparison theorem that does not use mixed weighted curvature. 
However, the curvature condition used in Theorem \ref{thm-LapComp-2} is the weighted flag curvature, which is much stronger than the mixed weighted Ricci curvature used in Theorem \ref{thm-LapComp-1}.

\begin{thm}\label{thm-LapComp-2}
	Let $(M,F,\mu)$ be a forward complete $n$-dimensional Finsler metric measure space with finite misalignment $\alpha$. Denote the forward distance function by $r$ and by $V$ a fixed vector field on $M$. Suppose 
	the weighted flag curvature $K^N$ of $M$ is bounded from below by $-\frac{N-1}{n-1}K$ with $K>0$, for some $N>n$. Then the nonlinear Laplacian of $r$ with reference vector $V$ satisfies 
	\begin{eqnarray}\label{Lap-r-6}
		\Delta^Vr\leq \alpha(N-1) \mathfrak{ct}_{-K}(r)+C_0,
	\end{eqnarray}
	where $C_0$ is a constant depending on $\alpha$, 
	and the bound of non-Riemannian curvatures $K_0$, namely, 
	$$F(U)+F^*(\mathcal{T})+F(\mathrm{div} C(V))\leq K_0,$$
	for any local vector fields $V,W$ in the definitions of $U$ and $\mathcal{T}$, which are given in \eqref{def-UyW} and \eqref{def-TVW}, respectively, and $\mathrm{div} C(V)=C^{ip}_{k\vert i}(V)V^k\frac{\delta}{\delta x^p}$.
\end{thm}

\begin{proof}
	We denote the forward distance function by $r$. It is known from the above equatily \eqref{Lap-r-1} that
	\begin{eqnarray}\label{Lap-r-2}
		\begin{split}
			\Delta^Vr&=g^{ij}(x,V)\Phi_i r_j+g^{ij}(x,V)\frac{\partial^2 r}{\partial x^i\partial x^j}+\frac{\partial}{\partial x^i}g^{ij}(x,V)r_j\\
			&=g^{ij}(x,V)\frac{\partial^2 r}{\partial x^i\partial x^j}+g^{ij}(x,V)\Phi_i r_j\\
			&\quad-[\Gamma^i_{im}(x,V)g^{mj}(x,V)+\Gamma^j_{im}(x,V)g^{mi}(x,V)+2C^{ij}_k(x,V)\nabla^V_iV^k]r_j,		
		\end{split}
	\end{eqnarray} 
	where $\nabla^V$ means the horizontal covariant derivative with respect to the Chern connection with the reference vector $V$. By the definition of the distortion, it follows from \eqref{Lap-r-2} that
	\begin{eqnarray}\label{Lap-r-3}
		\begin{split}
			\Delta^Vr=&\mathrm{tr}_{V}\nabla^{2}r+g^{ij}(x,V)r_k(\Gamma^k_{ij}(x,\nabla r)-\Gamma^k_{ij}(x,V))\\
			&-g^{ij}(x,V)\tau_i(x,V)r_j-2C^{ij}_k(x,V)\nabla^V_iV^kr_j,
		\end{split}
	\end{eqnarray}
	where $\mathrm{tr}_{V}\nabla^{2}r=g^{ij}(x,V)\nabla_jr_i$ means taking trace of the Hessian of $r$ with respect to the metric $g(x,V)$, and $\tau$ is the distortion. From the definition of the tensor $U_y(W)$, the second term on the RHS of \eqref{Lap-r-3} is $U_{V}(\nabla r)$. We denote the last term on the RHS of \eqref{Lap-r-3} by $2\mathrm{tr}_{\nabla V}C(V)(\nabla r)$, in which $\mathrm{tr}_{\nabla V}C(V)$ means taking trace of the Cartan tansor in the direction $V$ by $\nabla V$. Therefore, we get that 
	\begin{eqnarray}\label{Lap-r-4}
		\Delta^Vr=\mathrm{tr}_{V}\tilde{H}(r)+U_V(\nabla r)+g^*_V(\mathcal{T}(\nabla r, V),dr)-2\mathrm{tr}_{\nabla V}C(V)(\nabla r),
	\end{eqnarray}
	where $\tilde{H}(r)$ is the weighted Hessian of the forward distance $r$, whose definition is given in \eqref{Def-WHf}, $g^*_V$ is the dual metric with reference vector field $V$ and $\mathcal{T}(\nabla r, V)=\nabla^{\nabla r}\tau(\nabla r)-\nabla^V\tau(V)$. 
	
	
	The Cauchy-Schwarz inequality implies the following provided that the components of the tensor $U$ are $U^i=g^{kl}(V)(\Gamma^i_{kl}(x,\nabla r)-\Gamma^s_{kl}(x,V))$. 
	\begin{eqnarray}\label{est-U}
		\begin{split}
			U_V(\nabla r)&=g^{ij}(x,V)r_k(\Gamma^k_{ij}(x,\nabla r)-\Gamma^k_{ij}(x,V))\\
			&\leq [g_{ij}(\nabla r)U^i(V,\nabla r)U^j(V,\nabla r)]^{\frac12}[g^{ij}(\nabla r)r_ir_j]^{\frac12}\\
			&=F_{\nabla r}(U)\\
			&\leq \sqrt{\kappa} F(U),
		\end{split}
	\end{eqnarray}
	where $U=U^i\frac{\delta}{\delta x^i}$ is a vector field defined on the sphere bundle, and $\kappa$ is the smoothness constant. Also it can be deduced that
	\begin{eqnarray}\label{est-T}
		\begin{split}
			g_V^*(\mathcal{T}(\nabla r,V),dr)\leq& [g_{pq}(\nabla r)\mathcal{T}_i(V)g^{ip}(V)\mathcal{T}_j(V)g^{jq}(V)]^{\frac12}[g^{ij}(\nabla r)r_ir_j]^{\frac12}\\
			\leq&\sqrt{\alpha}[g_{pq}(V)\mathcal{T}_i(V)g^{ip}(V)\mathcal{T}_j(V)g^{jq}(V)]^{\frac12}\\
			=&\sqrt{\alpha}F^*(\mathcal{T}).
		\end{split}
	\end{eqnarray}
	And more,  
	\begin{eqnarray}\label{est-IVVr}
		\mathrm{tr}_{\nabla V}C(V)(\nabla r)&=-C^{ij}_{k\vert i}(V)V^kr_j\leq F_{\nabla r}(\mathrm{div} C(V))\leq \sqrt{\kappa} F(\mathrm{div}C(V)),
	\end{eqnarray}
	where $\mathrm{div} C(V)=C^{ip}_{\,\,\,\,k\vert i}(V)V^k\frac{\delta}{\delta x^p}$ is a vector field defined on the sphere bundle.
	
	Therefore, we can obtain from \eqref{Lap-r-4}-\eqref{est-IVVr} that
	\begin{eqnarray}\label{Lap-r-ineq}
		\Delta^Vr-\sqrt{\alpha}K_0\leq\mathrm{tr}_V\tilde{H}(r).
	\end{eqnarray}
	
	It follows from \eqref{Lap-r-ineq} and Theorem \ref{thm-HessComp} that
	\begin{eqnarray}
		\begin{split}
			\Delta^Vr-\sqrt{\alpha}K_0\leq&\sum_{i=1}^n\tilde{H}(r)(E_i,E_i)\\
			\leq&\frac{N-1}{n-1}\mathfrak{ct}_{-K}(r)\left(\sum_{i=1}^{n-1}g_{\nabla r}(E_i,E_i)-\sum_{i=1}^{n-1}g_{\nabla r}(\nabla r,E_i)^2\right)\\
			\leq&\frac{N-1}{n-1}\mathfrak{ct}_{-K}(r)\sum_{i=1}^{n-1}\frac{g_{\nabla r}(E_i,E_i)}{g_V(E_i,E_i)}\\
			\leq& \alpha(N-1)\mathfrak{ct}_{-K}(r).
		\end{split}
	\end{eqnarray}
\end{proof}



Now we turn to prove Theorem \ref{thm-LapComp-1}, starting from \eqref{Lap-r-3}.

\begin{proof}[Proof of Theorem \ref{thm-LapComp-1}.]
	

	To estimate the RHS of \eqref{Lap-r-ineq}, we first choose a $g_{V}$-orthonormal frame $\{e_1,\cdots,e_n\}$ at $x$. Without losing the generality, one can assume that $\nabla^{2}r$, the hessian of $r$, 
	is a diagonal matrix at $x$, by taking an orthogonal transformation of $\{e_i\}_1^n$ if necessary. $E_1,\cdots,E_{n-1}, E_n$ are obtained by parallel moving of $e_1,\cdots,e_{n-1}, e_n$ with respect to $g_{\nabla r}$ along the unique minimal geodesic $\gamma:[0,r(q)]\rightarrow M$ whence $r$ is smooth at $q\in M$. Thus we have
	\begin{eqnarray}\label{eq-drHEE}
		\begin{split}
			&\frac{d}{dr}(\nabla^2r(E_i,E_j))=\frac{d}{dr}g_{\nabla r}(\nabla_{E_i}\nabla r,E_j)\\
			=&-g_{\nabla r}(R_{\nabla r}(E_i,\nabla r)\nabla r,E_j)-g_{\nabla r}(\nabla_{\nabla _{E_i}\nabla r}\nabla r,E_j)\\
			=&-g_{\nabla r}(R_{\nabla r}(E_i,\nabla r)\nabla r,E_j)-\sum_kg_{V}(\nabla _{E_i}\nabla r,E_k)g_{\nabla r}(\nabla_{E^k}\nabla r,E_j)\\
			=&-g_{\nabla r}(R_{\nabla r}(E_i,\nabla r)\nabla r,E_j)-\sum_{k,l}g_V(E^k,E^l)g_{V}(\nabla _{E_i}\nabla r,E_k)g_{\nabla r}(\nabla_{E_j}\nabla r,E_l),
		\end{split}
	\end{eqnarray}
	where we denote $E^k=g_V^{-1}(E_k,E_l)E_l$ by taking $g_V^{-1}$ the inverse matrix of $g_V$ and use that 
	\begin{eqnarray}\label{eq-Eir}
		\nabla_{E_i}\nabla r=g_{V}(\nabla_{E_i}\nabla r,E_k)E^k,
	\end{eqnarray}
	in the third equality of \eqref{eq-drHEE},
	as well as the symmetry of the Hessian, that is, $g_{\nabla r}(\nabla_{E_k}\nabla r,E_j)=g_{\nabla r}(\nabla_{E_j}\nabla r,E_k)$.
	Taking the trace of \eqref{eq-drHEE} by $g_V$ yields
	\begin{eqnarray}\label{eq-drtrVr}
		\frac{d}{dr}(\mathrm{tr}_V\nabla^2r)=-\mathrm{tr}_{V}R_{\nabla r}(\nabla r)-\sum_{i,k}g_{V}(\nabla _{E_i}\nabla r,E^k)g_{\nabla r}(\nabla_{E_k}\nabla r,E^i).
	\end{eqnarray}
	
	Next, we claim that $\sum_{i,k}g_{V}(\nabla _{E_k}\nabla r,E^i)g_{\nabla r}(\nabla_{E_i}\nabla r,E^k)$ is positive. 
	By adopting the symmetry of $\nabla^2r$, we can directly get that
	\begin{eqnarray}
		\begin{split}
			&\sum_{i,k}g_{V}(\nabla _{E_k}\nabla r,E^i)g_{\nabla r}(\nabla_{E_i}\nabla r,E^k)
			\\
			=&\sum_{i,k,j}g_{\nabla r}(\nabla_{E_i}\nabla r,g_{V}(\nabla _{E_j}\nabla r,E_k)E_k)g_V(E^j,E^i)\\
			=&\sum_{i,j}g_{\nabla r}(\nabla_{E_i}\nabla r,\nabla_{E_j}\nabla r)g_V(E^j,E^i)\\
			=&:\|\nabla^2r\|_{HS(V,\nabla r)}^2,
		\end{split}
	\end{eqnarray}
	which we may call $\|\nabla^2r\|_{HS(V,\nabla r)}^2$ the \emph{mixed Hilbert-Schmit norm} with respect to $g_V$ and $g_{\nabla r}$. Moreover, if we choose $g_{V}$-orthonormal bases $\{\hat E_i\}$ and another $g_{\nabla r}$-orthonormal bases $\{\tilde E_i\}$ to express that 
	\begin{eqnarray}
		\|\nabla^2r\|_{HS(V,\nabla r)}^2=\sum_{i,j}g_{\nabla r}(\nabla_{\hat E_i}\nabla r,\tilde E_j)g_{\nabla r}(\nabla_{\hat E_i}\nabla r,\tilde E_j)=\sum_{i,j}(\nabla^2r(\hat E_i,\tilde E_j))^2,
	\end{eqnarray}
	which is obviously nonnegative.
	
	Denote the frames $\{\hat E_i\}$ and $\{\tilde E_i\}$ by $\hat E$ and $\tilde E$, respectively. There exists an invertible matrix $A$ such that $\tilde E=\hat EA$. By taking an orthogonal transformation, we may assume without loss of generality that $A^TA$ is a diagonal matrix with diagonal entries $0<\lambda_1\leq\lambda_2\leq\cdots\leq\lambda_n$. Denote the Hessian of $r$ by $G=\nabla^2r(\hat E,\hat E)$, a symmetric matrix. The Hilbert-Schmit norm of $G$ is $\|G\|^2_{HS(V)}=\sum_{i,j}G_{ij}G^T_{ji}=tr(GG^T)$, and the mixed Hilbert-Schmit norm of $\nabla^2r$ is
	$$\|\nabla^2r\|^2_{HS(V,\nabla r)}=\sum_{i,j,k}(\nabla^2r(\hat E_i,\hat E_k)A_{kj})^2=\sum_{i,j,k,l}G_{ik}A_{kj}A^T_{jl}G^T_{li}=tr(GAA^TG^T).$$
	From the choice of the basis $\hat E$ and $\tilde E$, we deduce that 
	\begin{eqnarray}\label{inequ-lambdaGG}
		\lambda_1 GG^T\leq GAA^TG^T\leq \lambda_n GG^T,
	\end{eqnarray}
	since the eigenvalues of $A^TA$ are the eigenvalues of $AA^T$. Equally,
	\begin{eqnarray}
		\lambda_1\|\nabla^2r\|^2_{HS(V)}\leq\|\nabla^2r\|^2_{HS(V,\nabla r)}\leq\lambda_n\|\nabla^2r\|^2_{HS(V)}.
	\end{eqnarray}
	\eqref{inequ-lambdaGG} means the inequality is valid for any vector $X$ that $$\lambda_1X^TGG^TX \leq X^TGAA^TG^TX \leq\lambda_nX^TGG^TX.$$
	Thus, more that the nonnegativity, we find that $\sum_{i,k} g_{V}(\nabla _{E_i}\nabla r,E^k) g_{\nabla r}(\nabla_{E_k}\nabla r,E^i)\geq \lambda_1\|\nabla^2r\|^2_{HS(V)}>0$. Hence, \eqref{eq-drHEE} provides that
	\begin{eqnarray}\label{ineq-drHEE}
		\frac{d}{dr}(tr_V\nabla^2r)\leq-tr_VR_{\nabla r}(\nabla r)-\lambda_1\|\nabla^2r\|_{HS(V)}^2.
	\end{eqnarray}
	
	Now we explain geometric meanings of $\lambda_1$ and $\lambda_n$. It follows from the transformation $\tilde E=\hat EA$ that
	\begin{eqnarray}\label{matrix-1}
		I=g_V(\hat E,\hat E)=g_V(\tilde EA^{-1},\tilde EA^{-1})=(A^T)^{-1}g_V(\tilde E,\tilde E)A^{-1},
	\end{eqnarray}
	where $I$ denotes the identity matrix. 
	Noticing the fact that $e=\tilde E^{-1}\frac{V}{F(V)}$ is an eigenvector of the eigenvalue $1$, one may deduce that $\lambda_1\leq1\leq\lambda_n$.
	In fact, \eqref{matrix-1} implies that
	\begin{eqnarray}
		A^TA=g_V(\tilde E,\tilde E).
	\end{eqnarray}
	Suppose $e$ is the eigenvector of $\lambda$, i.e., $e^Tg_{V}(\tilde E,\tilde E)e=\lambda=\lambda\vert e\vert ^2$, where $\vert \cdot\vert $ denotes the Euclidean norm. So we have 
	$$\lambda=\frac{e^Tg_{V}(\tilde E,\tilde E)e}{e^Tg_{\nabla r}(\tilde E,\tilde E) e}=\frac{g_{V}(\tilde Ee,\tilde Ee)}{g_{\nabla r}(\tilde Ee,\tilde Ee)}=\frac{g_V(\tilde e,\tilde e)}{g_{\nabla r}(\tilde e,\tilde e)},$$
	where $\tilde e=\tilde E e$.
	Then we know that
	\begin{eqnarray}
		\lambda_1=\inf_{g_{\nabla r}(e,e)=1}g_{V}(e,e) \quad\mbox{and}\quad \lambda_n=\sup_{g_{\nabla r}(e,e)=1}g_{V}(e,e).
	\end{eqnarray}
	
	According to the definition of the global misalignment with respect to $\nabla r$ (cf. Definition \ref{def-misalignment}), it represents that $\lambda_1=\alpha_m(x,\nabla r)$ and $\lambda_n=\alpha_M(x,\nabla r)$. Hence we have
	\begin{eqnarray}\label{ineq-alphalambda}
		\frac{1}{\alpha}\leq\lambda_1\leq1\leq\lambda_n\leq \alpha.
	\end{eqnarray}
	Therefore, it is deduced from \eqref{ineq-drHEE} and \eqref{ineq-alphalambda} that 
	\begin{eqnarray}\label{drtrVH}
		\frac{d}{dr}(\mathrm{tr}_V\nabla^2r)\leq-\mathrm{tr}_{V}R_{\nabla r}(\nabla r)-\frac{1}{\alpha}\|\nabla^2r\|^2_{HS(V)}.
	\end{eqnarray}
	
	It is easy to see that $\|\nabla^2r\|^2_{HS(V)}\geq \frac{(\mathrm{tr}_V\nabla^2r)^2}{n}$. 
	While one actually can get that $\|\nabla^2r\|^2_{HS(V)}\geq \frac{(\mathrm{tr}_V\nabla^2r)^2}{n-1}$ by noticing that one can choose a vector in direction $\nabla r$ with $\nabla^2r(E_n,E_n)=\nabla^2r(\nabla r,\nabla r)=0$.
	Combining with \eqref{Def-MWH} and the above fact, \eqref{drtrVH} means that
	\begin{eqnarray}\label{3.28}
		\begin{split}
			\frac{d}{dr}(\mathrm{tr}_V\tilde{H}(r))&\leq-\mathrm{tr}_{V}R_{\nabla r}(\nabla r)-\frac{d}{dr}S(\nabla r)-\frac{(\mathrm{tr}_V\tilde{H}(r)+S(\nabla r))^2}{\alpha(N-1)}\\
			&\leq-\mathrm{tr}_{V}R^{\nabla r}(\nabla r)-\frac{d}{dr}S(\nabla r)-\frac{(\mathrm{tr}_V\tilde{H}(r))^2}{N+(\alpha-1)n-\alpha}+\frac{(S(\nabla r))^2}{(N-n)}\\
			&=-\frac{(\mathrm{tr}_V\tilde{H}(r))^2}{N+(\alpha-1)n-\alpha}- {^mRic}^N_V(\nabla r),
		\end{split}
	\end{eqnarray}
	where $^mRic^N_V(\nabla r)=\mathrm{tr}_{V}R_{\nabla r}(\nabla r)+\dot S(\nabla r)-\frac{(S(\nabla r))^2}{(N-n)}$ is the mixed weighted Ricci curvature with tracing vector field $V$ and reference vector field $\nabla r$. We have adopted the inequality $(a\pm b)^2\geq\frac{a^2}{\lambda+1}-\frac{b^2}{\lambda}$, with $\lambda=\frac{N-n}{\alpha(N-1)}>0$. Denoting $C(N,\alpha)=N+(\alpha-1)n-\alpha$, one can estimate it as
	\begin{eqnarray}
		C(N,\alpha)=N-1+(\alpha-1)(n-1)\geq N-1.
	\end{eqnarray}
	
	Noticing the assumption of the mixed weighted Ricci curvature $^mRic^N_V(\nabla r)\geq -K$, and the fact that
	\begin{eqnarray}
		\frac{d}{dr}\mathfrak{ct}_{-l}(r)=-l-(\mathfrak{ct}_{-l}(r))^2,
	\end{eqnarray}
	with $l=\frac{K}{C(N,\alpha)}\leq\frac{K}{N-1}$. One can easily verify that
	\begin{eqnarray}
		\begin{split}
			\frac{d}{dr}\left[\mathrm{tr}_V\tilde{H}(r)-C(N,\alpha)\mathfrak{ct}_{-l}(r)\right]\leq -\frac{1}{C(N,\alpha)}\left[(\mathrm{tr}_V\tilde{H}(r))^2-(C(N,\alpha)\mathfrak{ct}_{-l}(r))^2\right].
		\end{split}
	\end{eqnarray}
	Canonically, we set 
	$$m(r)=(\mathrm{tr}_V\tilde{H}(r)-C(N,\alpha)\mathfrak{ct}_{-l}(r))e^{\frac{1}{C(N,\alpha)}\int_{\epsilon}^r(\mathrm{tr}_V\tilde{H}(r)+C(N,\alpha)\mathfrak{ct}_{-l}(r))d\tau},$$
	so that $\frac{d}{dr}m(r)\leq 0$. Therefore,
	\begin{eqnarray}\label{epsilon}
		\mathrm{tr}_V\tilde{H}(r)-C(N,\alpha)\mathfrak{ct}_{-l}(r)\leq\lim_{r\rightarrow 0+}(\mathrm{tr}_V\tilde{H}(r)-C(N,\alpha)\mathfrak{ct}_{-l}(r)).
	\end{eqnarray}
	Now we estimate the limitation on the RHS of \eqref{epsilon}. Let $p_0$ be the base point and $r(x)=d(p_0,x)$ be the forward distance from $p_0$ to $x$. Let $U$ be a neighborhood of $p_0$ such that $x\in U$ when $r$ is small enough. Suppose $a_1$ and $a_2$ are the lower and upper bounds of the flag curvature, respectively. The Hessian of the distance function $r$ satisfies that 
	\begin{eqnarray}
		\mathfrak{ct}_{a_2}(r)\leq\nabla^2r(X,X)\leq \mathfrak{ct}_{a_1}(r),
	\end{eqnarray}
	for any vector $X$. From the definition of the mixed weighted Hessian \eqref{Def-MWH}, one deduces that 
	\begin{eqnarray}
		\tilde{H}^V(r)(X,X)=\nabla^2r(X,X)-\frac{S(\nabla r)}{n-1},
	\end{eqnarray}
	for any $X$ satisfying $g_V(X,\nabla r)=0$ and $g_V(X,X)=1$. Therefore, taking $\{X_i\}$ as an orthonormal basis with respect to $g_V$ and orthogonal to $\nabla r/\sqrt{g_V(\nabla r,\nabla r)}$,
	\begin{eqnarray}
		\mathrm{tr}_V\tilde{H}(r)=\sum_{i=1}^{n-1}\tilde{H}^V(r)(X_i,X_i)\leq\frac{n-1}{r}-\frac{n-1}{3}a_1r-S(\nabla r),
	\end{eqnarray}
	which implies that 
	\begin{eqnarray*}
		\mathrm{tr}_V\tilde{H}(r)-C(N,\alpha)\mathfrak{ct}_{-l}(r)\leq\frac{n-N}{r}+\frac13\left(\frac{N-1}{C(N,\alpha)}K-(n-1)a_1\right)r-S(\nabla r).
	\end{eqnarray*}
	So there exists $\epsilon>0$ such that $\mathrm{tr}_V\tilde{H}(r)-C(N,\alpha)\mathfrak{ct}_{-l}(r)<0$, when $r<\epsilon$ and $n<N$. Hence we get 
	\begin{eqnarray}\label{ineq-trVr}
		\mathrm{tr}_V\tilde{H}(r)\leq C(N,\alpha)\mathfrak{ct}_{-l}(r),
	\end{eqnarray}
	for any $r>0$. Combining \eqref{Lap-r-ineq} and \eqref{ineq-trVr}, we get that
	\begin{eqnarray}
		\Delta^Vr\leq C(N,\alpha)\mathfrak{ct}_{-l}(r)+C_0,
	\end{eqnarray}
	by taking $C_0=\sqrt{\alpha}K_0$.
\end{proof}


Remark \ref{rem-localalpha} suggests that when applying this comparison theorem, we can relax the global bounds of misalignment $\alpha$ to local bounds in order to include a larger range of Finsler manifolds.\\

\begin{rem}\label{rem-comparison}
	\begin{itemize}
		\item[i)] 	According to \eqref{3.28}, one may find that the mixed weighted Ricci curvature $^mRic^N_V(\nabla r)\geq -K(2R)$ can be substituted by $^mRic^{\infty}_V(\nabla r)\geq -K(2R)$ and the $S$-curvature bounded, i.e., $\vert S\vert \leq K'(2R)$ for some constant $K'$ on $B(2R)$.
		\item[ii)] We will not use this comparison theorem until Section \ref{Ch-noncpt}.\\
	\end{itemize} 
\end{rem}

It is an obvious fact that the Mixed weighted Ricci curvature condition is stronger than the weighted Ricci curvature condition. A natural and important question arises here, that is, what is the real relationship between the Mixed curvature condition and the curvature dimension condition? A detailed study of this issue will be further examined in another article.

\section{Finslerian Schr\"odinger equation}\label{SecFSE}

The main task of this section is to introduce the Finslerian Schr\"odinger equation. Before that, let's briefly review some of the research history of Schr\"odinger's equations related to our topic.

Schr\"odinger equation has important applications in both mathematical and physical theory. Since Li-Yau's work, there has been a lot of work on the study of Schr\"odinger operators on manifolds.
R. Chen obtained the lower bound of the eigenvalues of the elliptic Schr\"odinger operator of the Neumann boundary condition and the lower bound of the gap between eigenvalues in the case of Ricci's curvature and the second fundamental form having a lower bound \cite{Chen1997PJM}. T. H. Colding and W. P. Minicozzi II discussed the asymptotic behavior of the solution to the elliptic Schr\"odinger equation on Riemannian manifolds with non-negative Ricci curvature and Euclidean volume growth, as well as its applications \cite{CM1997AJM}. A. Grigor'yan and W. Hansen studied the Liouville property of the elliptical Schr\"odinger equation on Riemannian manifolds \cite{GH1998MA}. Q. S. Zhang gave estimates of the global upper and lower bounds of the parabolic Schr\"odinger kernels on the Riemannian manifold with non-negative Ricci curvature \cite{QSZhang2000CMP}. He then proved certain existential conditions for elliptical Schr\"odinger heat kernels and obtained certain comparative properties \cite{QSZhang2003BLMS}. Later, together with P. Souplet, he also obtained a gradient estimate of the heat equation on noncompact manifolds \cite{SZ2006BLMS}. Inspired by this, many researchers have made significant improvements  \cite{DK2015AM,Z2016AMS,CZ2018JMAA,MZ2018CRMASP}. In addition, E. M. Ouhabaz and C. Poupaud proved the Cwikel-Lieb-Rozenblum and Lieb-Thirring type estimates of the elliptic Schr\"odinger operator on complete Riemannian manifolds \cite{OP2010AAM}. Recently, F. Faraci and C. Farkas obtained a bounded characteristic of Schr\"odinger equations on complete noncompact Riemannian manifolds \cite{FF2019CCM}. L. Zhang studied matrix and classical Harnack estimates for positive solutions to the nonlinear diffusion equation on a compact K\"ahler manifold with fixed metric or under the normalized K\"ahler-Ricci flow \cite{ZhangLD2020}. J. Wu  gave new Gaussian type upper bounds for the Schr\"odinger heat kernel on complete gradient shrinking Ricci solitons with the scalar curvature bounded above \cite{Wu2021JGP}. X. Huang, Y. Sire and C. Zhang studied the eigenfunction $L^p$ bounds estimation of elliptic spectral fractional Schr\"odinger operators on closed Riemannian manifolds, and obtained Strichartz estimates for the fractional wave equation \cite{HSZ2021JMPA}.  C. Cavaterra, S. Dipierro, Z. Gao and E. Valdinoci obtained a broader global gradient estimate of a class of nonlinear parabolic equations \cite{CDGV2022JGA}. K. Merz considered fractional Schr\"odinger operators with possibly singular potentials and derive certain spatially averaged estimates for its complex-time heat kernel \cite{Merz2022}. X. Huang and C. Zhang also studied the pointwise Weyl Law for the Schr\"odinger operators on compact Riemannian manifolds without boundary \cite{HZ2022AM}.

We consider the \emph{Finslerian Schr\"odinger equation} \eqref{equ-FS1}
\begin{eqnarray*}
	(\Delta^{\nabla u}-\partial_t-q)u=0,
\end{eqnarray*}
on a Finsler metric measure space $(M,F,\mu)$, where $\Delta^{\nabla u}$ is the nonlinear Laplacian given in \eqref{def-WLap} and $q=q(x,t)\in C([0,T],H^1(M))$. The solution to \eqref{equ-FS1} is a function $u(x)$ on $M$. The regularity of function $q(x,t)$ is required to be continuous in $t$ for the solution existence, while we only need $q(x,t)\in L^2([0,T],H^1(M))$ for the gradient estimation of the positive solutions.

We also use $\nabla$ to denote the gradient and $\Delta$ to denote the Laplacian when they act on function $h$ in direction $\nabla h$. 
One may notice that $\nabla$ also denotes 
the horizontal Chern connection with reference vector $\nabla h$. 
However, for any function $u$ on manifold $M$, the two concepts are the same, i.e.,
$$\nabla_i u(x)=\frac{\delta u(x)}{\delta x^i}=\frac{\partial u}{\partial x^i}=u_i=g_{\nabla u}(\nabla u,\frac{\partial}{\partial x^i}),$$ 
where $\nabla_i u$ means the horizontal Chern derivative locally along $x^i$ and $\nabla u$ means the gradient of $u$. 

Furthermore, the Hessian of $u$ is defined by 
$$\nabla^2u(X,Y)=g_{\nabla u}(D^{\nabla u}_X\nabla u,Y),$$
for any $X,Y\in TM$, where $D^{\nabla u}$ is the covariant differentiation with respect to the horizontal Chern connection, whose reference vector is $\nabla u$. By the definition, it is easy to check that $\nabla^2u$ is a symmetric 2-form.

Let $H^1(M):=W^{1,2}(M)$ and $H^1_0(M)$ is the closure of $C^{\infty}_0(M)$ under the norm 
$$\|u\|_{H^1}:=\|u\|_{L^2(M)}+\frac12\|F(\nabla u)\|_{L^2(M)}+\frac12\|\overleftarrow{F}(\overleftarrow{\nabla} u)\|_{L^2(M)},$$
where $\overleftarrow{F}$ is the \emph{reverse Finsler metric}, defined by $\overleftarrow{F}(x,y):=F(x,-y)$ for all $(x,y)\in TM$, and $\overleftarrow{\nabla} u$ is the gradient of $u$ with respect to the reverse metric $\overleftarrow{F}$. In fact, $\overleftarrow{F}(\overleftarrow{\nabla} u)=F(\nabla(-u))$. Then $H^1(M)$ is a Banach space with respect to the norm $\|\cdot\|_{H^1}$.

Analog to \cite{Ohta2009-2}, one can define that
\begin{defn}
	For $T>0$, a function $u$ on $[0,T]\times M$ is a global solution to the Finslerian Schr\"odinger equation \eqref{equ-FS1} if $u\in L^2([0,T], H^1_0(M))\cap H^1([0,T],H^{-1}(M))$, and if, for every time $t\in [0,T]$, any test function $\varphi\in H^1_0(M)$ (or $\varphi\in C^{\infty}_0(M)$), it holds that
	$$\int_M\varphi(u_t+qu)d\mu=-\int_{M}d\varphi(\nabla u)d\mu.$$
\end{defn}

The proof of the following interior regularity of the solution on $[0,\infty)\times M$ is the same as the one in \cite{Ohta2009-2}.

\begin{thm}
	Suppose $(M,F,\mu)$ is a Finsler metric measure space with finite reversibility $\rho<\infty$. Then one can take the continuous version of a global solution $u$ of the Finslerian Schr\"odinger equation \eqref{equ-FS1}, and it enjoys the $H_{loc}^2$-regularity in $x$ as well as the $C^{1,\beta}$-regularity in both $t$ and $x$ for some $0<\beta<1$. Moreover, $u_t$ lies in $H^1_{loc}(M)\cap C(M)$, and further in $H^1_0(M)$ if $F$ has finite uniform smoothness constant $\kappa$. The elliptic regularity shows that $u$ is $C^{\infty}$ on $\cup_{t>0}(\{t\}\times M_{u(t,x)})$.	
\end{thm}

Given an open subset $\Omega\subset M$ and an open interval $I\subset \mathbb{R}$, a local solution to the Finslerian Schr\"odinger equation \eqref{equ-FS1} is defined by
\begin{defn}\label{def-localsolution}
	A function $u$ on any open subset $\Omega\subset M$ is a local solution to the Finslerian Schr\"odinger equation \eqref{equ-FS1}, if $u\in L^2_{loc}(I\times \Omega)$ with $F^*(du)\in L^2_{loc}(I\times \Omega)$, and if, for any local test function $\varphi\in H^1_0(I\times \Omega)$ (or $\varphi\in C^{\infty}_0(I\times \Omega)$), it holds that
	$$\int_I\int_{\Omega}(u\varphi_t-qu\varphi) d\mu dt=\int_I\int_{\Omega}d\varphi(\nabla u)d\mu dt.$$\\
\end{defn}

Let $B_{2R}:=B^+_{2R}(x_0)$ be a forward geodesic ball with radius $2R$ in $M$ centered at $x_0$, and $I$ be an open interval in $\mathbb{R}$. As the global solution, if $u$ is a local solution to the Schr\"odinger equation \eqref{equ-FS1} on $I\times B_{2R}$, then $u\in H^2(B_{2R})\cap C^{1,\beta}(I\times B_{2R})$ 
with $\Delta u\in H^1(B_{2R})\cap C(B_{2R})$. Moreover $u_t$ lies in $H^1(B_{2R})\cap C^{\beta}(B_{2R})$, and further in $H^1_0(M)$ if $F$ has finite uniform smoothness constant.

\begin{rem}
	\begin{itemize}
		\item[(i)] $u$ is a solution to the Finslerian Schr\"odinger equation $(\Delta-\partial_t-q)u=0$ in the distributional sense if and only if $-u$ is a solution to $(\overleftarrow{\Delta}-\partial_t-q)u=0$ in the distributional sense, where $\overleftarrow{\Delta}$ is the Finsler Laplacian associated to $\overleftarrow{F}$.
		Since the metric $\overleftarrow{F}$ is also a Finsler metric on $M$, so the results we obtained in this manuscript also satisfy for the solutions to $(\overleftarrow{\Delta}-\partial_t-q)u=0$. More precisely, we can use positive solutions to $(\overleftarrow{\Delta}-\partial_t-q)u=0$ to obtain the result for bounded negative solutions to $(\Delta-\partial_t-q)u=0$.
		\item[(ii)] The test function $\varphi$ can be chosen only vanishing on the spacial boundary, namely, $\varphi\in H^1(I\times B_R)$ with $\varphi(t,\cdot)\in H^1_0(B_R)$ for all $t\in I$.\\
	\end{itemize}
\end{rem}

The regularity of the time variable $t$ of the solution $u$ can be improved by the regularization $u_{\epsilon}$, which is defined by 
$$u_{\epsilon}(t,x)=J_{\epsilon}*u(t,x)=\int_{I}J_{\epsilon}(t-s)u(s,x)ds,\quad t\in I_{\epsilon}=\{t\in I \mid \mathrm{dist}(t,\partial I)>\epsilon\}$$
for any $\epsilon>0$ and $x\in B_{2R}$, where $J_{\epsilon}$ is defined by
\begin{eqnarray}
	J_{\epsilon}(t)=\frac{1}{\epsilon}\begin{cases}
		& k\exp\left(\frac{\epsilon^2}{t^2-\epsilon^2}\right), \quad \mbox{if } \vert t\vert <\epsilon,\\
		&0,\quad\quad\quad\quad\quad\quad\,\, \mbox{if }\vert t\vert \geq\epsilon,
	\end{cases}
\end{eqnarray}
where $k>0$ is chosen such that $\int_{\mathbb{R}}J_{\epsilon}(t)dt=1$. Then $u_{\epsilon}(t,x)$ is smooth in $t\in I_{\epsilon}$, and $\lim_{\epsilon\rightarrow 0}u_{\epsilon}(t,x)=u(t,x)$ for a.e. $t\in I_{\epsilon}$ and uniformly on each compact subsets of $I$ for each $x\in B_{2R}$. In particular,we can take $I_{\epsilon}=(\epsilon,T-\epsilon)$ when $I=(0,T)$ for any $0<T<\infty$. 

The following lemma is an analog of the one in \cite{Xia2023}.
\begin{lem}\label{lem-Ju}
	Let $u$ be a local solution to the Finslerian Schr\"odinger equation \eqref{equ-FS1} on $I\times B_{2R}$. Then for any $\epsilon>0$, we have
	\begin{itemize}
		\item[(i)] $\lim_{\epsilon\rightarrow0}(u_{\epsilon})_t=u_t$, $\lim_{\epsilon\rightarrow0}\nabla u_{\epsilon}=\nabla u$ and $\lim_{\epsilon\rightarrow 0}F(\nabla u_{\epsilon})=F(\nabla u)$ on $I_{\epsilon}\times B_{2R}$;
		\item[(ii)] If $u\in L^2(I\times B_{2R})$ with $F^*(du)\in L^2(I\times B_{2R})$, then $u_{\epsilon}\in L^2(I_{\epsilon}\times B_{2R})$ with $F^*(du_{\epsilon})\in L^2(I_{\epsilon\times B_{2R}})$. Further, $(u_{\epsilon})_t\in L^2(I_{\epsilon}\times B_{2R})$ and hence $u_{\epsilon}\in H^1(I_{\epsilon}\times B_{2R})\cap C^{\infty}(\mathbb{R})\cap H^2(B_{2R})$;
		\item[(iii)] $u_{\epsilon}$ is a local solution to the Schr\"odinger equation \eqref{equ-FS1} on $I_{\epsilon}\times B_{2R}$.
	\end{itemize}
\end{lem}

\begin{proof}
	We only need to prove (iii). The proofs of (i) and (ii) are the same to Lemma 4.1 in \cite{Xia2023}. Since $J_{\epsilon}$ is a smooth function with compact support and satisfies $\partial_t(J_{\epsilon}(t-s))=-\partial_{s}(J_{\epsilon}(t-s))$, one can get the following according to the Fubini's Theorem for any $\phi(t,x)\in H_0^1(I_{\epsilon}\times B_{2R})$.
	\begin{eqnarray}
		\begin{split}
			\int_{I_{\epsilon}}\int_{B_{2R}}\phi(t,x)\partial_{t}u_{\epsilon}dmdt&=\int_{I_{\epsilon}}\int_I\int_{B_{2R}}\phi(t,x)\partial_t(J_{\epsilon}(t-s))u(s,x)dmdsdt\\
			&=-\int_{I_{\epsilon}}\int_I\int_{B_{2R}}\partial_s(\phi(t,x)J_{\epsilon}(t-s))u(s,x)dmdsdt\\
			&=-\int_{I_{\epsilon}}\int_I\int_{B_{2R}}[J_{\epsilon}(t-s)d\phi(\nabla u(s,x))\\
			&\quad\quad\quad\quad\quad+J_{\epsilon}(t-s)\phi(t,x)q(s,x)u(s,x)]dmdsdt\\
			&=-\int_{I_{\epsilon}}\int_{B_{2R}}[\phi(t,x)q(t,x)u_{\epsilon}(t,x)+d\phi(\nabla u_{\epsilon}(t,x))]dmd,
		\end{split}
	\end{eqnarray}
	where we have used Definition \ref{def-localsolution} with the test function $\phi(t,x)J_{\epsilon}(t-s)$ in the third equality.
\end{proof}

When studying the well-posedness of solutions to such nonlinear Schr\"odinger equations, in addition to the regularity of the solutions briefly described in this section, two natural and important topics are the existence of the solutions of the equations, as well as the uniqueness of the solution. According to the quasilinear of the equation \eqref{equ-FS1}, these properties of the solution can be demonstrated by the fixed point method. Issues related to these problems are also deviated from the keynote of our topic, which we will address in a dedicated article.

\section{Global gradient estimates on compact Finsler manifolds}\label{Ch-GGE}

In this section, we will show the global gradient estimates of positive solutions to \eqref{equ-FS1} on compact Finsler manifolds, with two different  weighted Ricci curvatures bounded from below. The curvature condition is first given by Ohta in \cite{Ohta2009-2} and \cite{Ohta2014} to generalize the Li-Yau gradient estimates of the positive solutions to heat equation on compact Finsler manifolds. Later, it has been widely adopted by researchers to study the Finsler geometric analysis problems.

\subsection{Weighted Ricci curvature $Ric^{N}$ bounded from below}\label{Ch-GGE-1}

Let's denote $u$ as a positive solution on $M\times [0,T]$ to the Finslerian Schr\"odinger equation \eqref{equ-FS1}. The function $u$ may not have enough regularity in both $t$ and $x$. However, Lemma \ref{lem-Ju} suggest us to utilize the smooth approximation $u_{\epsilon}$ in $t$ to overcome this problem. Thus, we can consider the equation on $M_u\times [0,T]$, where $M_u=\{x\in M\mid \nabla u(x)\neq 0\}$. It satisfies the following lemma.
\begin{lem}\label{lem-f}
	Let $u$ be a positive solution on $M\times [0,T]$ to the Finslerian Schr\"odinger equation \eqref{equ-FS1}, and $f(x,t)=\log u$. $f$ satisfies that 
	\begin{eqnarray}\label{equ-f-s}
		(\partial_t-\Delta)f=F^2(\nabla f)-q,
	\end{eqnarray}
	on $M_u\times [0,T]$.
\end{lem}
\begin{proof}
	Direct computation gives that
	\begin{eqnarray*}
		\mathrm{tr}_{\nabla u}\nabla^2f=\frac{\mathrm{tr}_{\nabla u}\nabla^2u}{u}-\frac{F^2(\nabla u)}{u^2}.
	\end{eqnarray*}
	Therefore,
	\begin{eqnarray*}
		\begin{split}
			\Delta f&=\mathrm{tr}_{\nabla u}\nabla^2f-S(\nabla f)\\
			&=\frac{\mathrm{tr}_{\nabla u}\nabla^2u}{u}-\frac{F^2(\nabla u)}{u^2}-S(\frac{\nabla u}{u})\\
			&=\frac{\Delta u}{u}-F^2(\nabla f),
		\end{split}
	\end{eqnarray*}
	in the distributional sense on $M_u$. We can get \eqref{equ-f-s} by noticing \eqref{equ-FS1}. 
\end{proof}
Another form of \eqref{equ-f-s} is the following.
\begin{cor}
	Let $u$ be a positive solution on $M\times [0,T]$ to \eqref{equ-FS1}, then $f(x,t)=\log u$ satisfies that
	\begin{eqnarray}\label{equ-f}
		\int_0^T\int_M[d\varphi(\nabla f)-\varphi_tf]d\mu dt=\int_0^T\int_M\varphi(F^2(\nabla f)-q)d\mu dt,
	\end{eqnarray}
	for any test function $\varphi\in C^2(M\times [0,T])$.
\end{cor}

Before the proof of the estimate, we introduce the following operators commutative lemma.
\begin{lem}\label{lem-com-tNtD}
	The nonlinear First- and second-order operators $\nabla^{\nabla u}$ and $\Delta^{\nabla u}$ is commutative with $\partial_t$. Namely, we have 
	\begin{eqnarray}
		[\partial_t,\nabla^{\nabla u}]\xi(u)=[\partial_t,\Delta^{\nabla u}]\xi(u)=0,
	\end{eqnarray}
	for any $C^1$ function $\xi$.
\end{lem}
\begin{proof}
	The proof is just some direct computations including
	\begin{eqnarray}
		\begin{split}
			\frac{\partial}{\partial t}\nabla^{\nabla u}\xi(u)=&\xi'\frac{\partial}{\partial t}\nabla u\\
			=&\xi'\frac{\partial}{\partial t}(g^{ij}(\nabla u)u_i\frac{\partial}{\partial x^j})\\
			=&\xi'[-2C^{ijk}(\nabla u)(u_k)_tu_i\frac{\partial}{\partial x^j}+g^{ij}(\nabla u)(u_i)_t]\frac{\partial}{\partial x^j}\\
			=&[g^{ij}(\nabla u)\xi'(u_t)_i]\frac{\partial}{\partial x^j}=\nabla^{\nabla u}(\frac{\partial}{\partial t}\xi(u)),
		\end{split}
	\end{eqnarray}
	and 
	\begin{eqnarray}
		\begin{split}
			\frac{\partial}{\partial t}\Delta^{\nabla u}\xi(u)=&\frac{\partial}{\partial t}[e^{-\Phi}\frac{\partial}{\partial x^i}(\xi'u^ie^{\Phi})]\\
			=&e^{-\Phi}\frac{\partial}{\partial x^i}(\xi'\frac{\partial}{\partial t}u^ie^{\Phi})\\
			=&e^{-\Phi}\frac{\partial}{\partial x^i}\{\xi'[-2C^{ijk}(\nabla u)(\frac{\partial}{\partial t}u_k)u_j+g^{ij}(\nabla u)(\frac{\partial}{\partial t}u_j)]e^{\Phi}\}\\
			=&e^{-\Phi}\frac{\partial}{\partial x^i}(g^{ij}(\nabla u)\xi'\frac{\partial}{\partial x^i}u_te^{\Phi})=\Delta^{\nabla u}\frac{\partial}{\partial t}(\xi(u)).
		\end{split}
	\end{eqnarray}
\end{proof}	

Let $H=t(F^2(\nabla f)-\beta f_t-\beta q)$ be a new function. Then $H$ satisfies the following lemma according to \eqref{equ-f-s}.
\begin{lem}\label{lem-H}
	Let $(M,F,\mu)$ be a Finsler metric measure space with weighted Ricci curvature $Ric^N$ having lower bound $-K$, where $K\geq0$ is a nonnegative constant. For any $\beta\geq 1$, $H=t(F^2(\nabla f)-\beta f_t-\beta q)$ satisfies
	\begin{eqnarray}
		\begin{split}
			H_t-\Delta^{\nabla u}H\leq 2[&df(\nabla^{\nabla u}H)+\frac{H}{2t}+KtF^2(\nabla f)-\frac{t}{N}(F^2(\nabla f)-f_t-q)^2\\
			&+(\beta-1)tdf(\nabla^{\nabla u}q)]-\beta t\Delta^{\nabla u}q,
		\end{split}
	\end{eqnarray}
	on $M_u\times (0,T]$.
\end{lem}
\begin{proof}
	We present the computations on $M_u\times(0,T]$.
	From the definition of $H$, we see that
	\begin{eqnarray}
		\nabla^{\nabla u} H=t[\nabla^{\nabla u}(F^2(\nabla f))-\beta \nabla^{\nabla u} f_t-\beta\nabla^{\nabla u} q].
	\end{eqnarray}
	Moreover, it satisfies that
	\begin{eqnarray*}
		\begin{split}
			\mathrm{tr}_{\nabla u}\nabla^{\nabla u 2}H&=t[2df(\mathrm{tr}_{\nabla u}\nabla^2(\nabla f))+2\|\nabla^2f\|_{HS(\nabla u)}-\beta \mathrm{tr}_{\nabla u}\nabla^2 f_t-\beta \mathrm{tr}_{\nabla u}\nabla^{\nabla u 2}q]\\
			&=t[2df(\nabla^{\nabla u}(\mathrm{tr}_{\nabla u}\nabla^2f))+2Ric(\nabla f,\nabla f)+2\|\nabla^2f\|_{HS(\nabla u)}\\
			&\quad\quad-\beta \mathrm{tr}_{\nabla u}\nabla^2 f_t-\beta \mathrm{tr}_{\nabla u}\nabla^{\nabla u 2}q],
		\end{split}
	\end{eqnarray*}
	on $M_u\times(0,T]$. By a direct computation in Finsler geometry, one may find that
	\begin{eqnarray}\label{Delta H-1}
		\begin{split}
			\Delta^{\nabla u}H&=\mathrm{div}_{\mu}\nabla^{\nabla u} H=e^{-\Phi}\frac{\partial}{\partial x^k}\left(e^{\Phi}g^{kl}(\nabla u)\frac{\partial H}{\partial x^l}\right)\\
			&=g_{\nabla u}(\nabla^{\nabla u}\Phi,\nabla^{\nabla u} H)+\mathrm{tr}_{\nabla u}\nabla^{\nabla u 2}H-\frac12g_{\nabla u}(\nabla^{\nabla u}\log\det(g_{ij}),\nabla^{\nabla u} H)\\
			&\quad+2C_{\nabla^2u}^{\nabla u}(\nabla^{\nabla u}H)+2P_{\nabla^2 u}(\nabla f,\nabla f),
		\end{split}
	\end{eqnarray}
	where it has utilized the fact that $C(\nabla f,\cdot,\cdot)=0$ and $C_{\nabla^2u}^{\nabla u}(\nabla^{\nabla u}H)=u^m_{\,\,\,\,\,\vert k}C^{kl}_{m}(\nabla u)H_l$ and $P_{\nabla^2 u}(\nabla f,\nabla f)=g^{ij}u^k_{\,\,\,\vert i}P^{\,\,\,\,\,l}_{m\,\,jk}f_lf^m$. The notation $``\vert "$ means the horizontal derivative with respect to the Chern connection in the direction $\nabla u$, and the directions in Cartan tensor and Landsberg curvature are both $\nabla u$.
	However, noticing that we always taking the direction in the non-Riemannian tensors to be $\nabla u$, thus 
	$$P_{\nabla^2 u}(\nabla f,\nabla f)=-\frac{1}{u}u^k_{\,\,\,\vert j}L^{ij}_kf^k=0,$$
	whence $\nabla u\neq 0$.
	Therefore, \eqref{Delta H-1} means that
	\begin{eqnarray}\label{Delta H-2}
		\Delta^{\nabla u}H=-d\tau(\nabla^{\nabla u} H)+\mathrm{tr}_{\nabla u}\nabla^{\nabla u 2}H+2C_{\nabla^2u}^{\nabla u}(\nabla^{\nabla u}H)
	\end{eqnarray}
	is valid on $M_u$, where $\tau$ is the distortion of the Finsler metric measure space $(M,F,\mu)$. Similarly, we also obtain on $M_u$ that
	\begin{align*}
		\Delta f&=-d\tau(\nabla f)+\mathrm{tr}_{\nabla u}\nabla^2f,\\
		\Delta^{\nabla u}q&=-d\tau(\nabla^{\nabla u} q)+\mathrm{tr}_{\nabla u}\nabla^{\nabla u 2}q+2C_{\nabla^2u}^{\nabla u}(\nabla^{\nabla u} q),
	\end{align*}
	where $C_{\nabla^2u}^{\nabla u}(\nabla^{\nabla u} q)=u^m_{\,\,\,\,\,\vert k}C^{kl}_{m}(\nabla u)q_l$. Combining these equalities with \eqref{Delta H-2} yields
	\begin{eqnarray}\label{Delta H-3}
		\begin{split}
			\Delta^{\nabla u}H-2C_{\nabla^2u}^{\nabla u}(\nabla^{\nabla u}H)=t\{&2df(\nabla^{\nabla u} \Delta f)+2g_{\nabla u}(\nabla^{\nabla u2}\tau,\nabla f\otimes\nabla f)\\
			&+2Ric(\nabla f,\nabla f)+2\|\nabla^{2}f\|^2_{HS(\nabla u)}-\beta(\Delta f)_t\\
			&-\beta(\Delta^{\nabla u}q-2C_{\nabla^2u}^{\nabla u}(\nabla^{\nabla u} q))\}\\
			=t\{&2df(\nabla^{\nabla u} [f_t-F^2(\nabla f)+q])+2\|\nabla^{2}f\|^2_{HS(\nabla u)}\\
			&+2Ric(\nabla f,\nabla f)+{\nabla u}(\nabla^{\nabla u2}\tau,\nabla f\otimes\nabla f)\\
			&-\beta(f_t-F^2(\nabla f)+q)_t\\
			&-\beta(\Delta^{\nabla u}q-2C_{\nabla^2u}^{\nabla u}(\nabla^{\nabla u} q))\},
		\end{split}
	\end{eqnarray}
	on $M_u$ according to Lemma \ref{lem-f}. The following calculations are straightforward and accessible to anyone familiar with Finsler geometry.
	\begin{eqnarray}\label{trNCH}
		\begin{split}
			C_{\nabla^2u}^{\nabla u}(\nabla^{\nabla u}H)&=u^m_{\,\,\,\,\,\vert l}C^l_{mk}t(2f^qf_q^{\,\,\,\vert k}-\beta f_t^{\,\,\,\vert k}-\beta q^k)\\
			&=2tu^m_{\,\,\,\,\,\vert l}[(C^l_{mk}f^k)_{\vert q}-C^l_{mk\vert q}f^k]f^q-\beta u^m_{\,\,\,\,\,\vert l}C^l_{mk}t(f^k_{\,\,\,t}+q^k)\\
			&=2tu^m_{\,\,\,\,\,\vert l}L^l_{mk}f^k-\beta u^m_{\,\,\,\,\,\vert l}t[(C^l_{mk}f^k)_t+C^l_{mk}q^k]\\
			&=-2t\beta C_{\nabla^2u}^{\nabla u}(\nabla^{\nabla u} q),
		\end{split}
	\end{eqnarray}
	where $``\vert "$ means the horizontal derivative with respect to the Chern connection in the direction $\nabla u$. 
	Thus, with the assistance of \eqref{trNCH}, as well as the definition of the S-curvature and the weighted Ricci curvature, \eqref{Delta H-3} is equal to
	\begin{eqnarray}\label{Delta h-4}
		\begin{split}
			\Delta^{\nabla u}H=t\{&2df(\nabla^{\nabla u}[f_t-F^2(\nabla f)+q])+2Ric^{\infty}(\nabla f,\nabla f)\\
			&+2\|\nabla^2f\|^2_{HS(\nabla u)}-\beta[(f_t-F^2(\nabla f)+q)_t+\Delta^{\nabla u}q]\},
		\end{split}
	\end{eqnarray}
	Utilizing the equality $\frac{(a+b)^2}{n}=\frac{a^2}{N}-\frac{b^2}{N-n}+\frac{N(N-n)}{n}(\frac{a}{N}+\frac{b}{N-n})^2$, for any $N>n$, one could obtain that
	\begin{eqnarray*}
		\|\nabla^2f\|^2_{HS(\nabla u)}=\frac{(\Delta f+ S(\nabla u))^2}{n}\geq \frac{(\Delta f)^2}{N}-\frac{( S(\nabla f))^2}{N-n},
	\end{eqnarray*}
	for any $N>n$. Therefore, \eqref{Delta h-4} implies that 
	\begin{eqnarray}\label{Delta h-5}
		\begin{split}
			\Delta^{\nabla u}H\geq& t\left\{2df(\nabla^{\nabla u} \Delta f)+2Ric^{N}(\nabla f,\nabla f)+2\frac{(\Delta f)^2}{N}-\beta[(\Delta f)_t+\Delta^{\nabla u}q]\right\}\\
			=&t\{2df(\nabla^{\nabla u} [f_t-F^2(\nabla f)+q])+2Ric^{N}(\nabla f,\nabla f)\\
			&\quad+2\frac{(f_t-F^2(\nabla f)+q)^2}{N}-\beta[(f_t-F^2(\nabla f)+q)_t+\Delta^{\nabla u}q]\}
		\end{split}
	\end{eqnarray}
	on $M_u$. Utilizing Lemma \ref{lem-f}, as well as the definition of $H$, we see 
	$\Delta f=-F^2(\nabla f)+q+f_t=-\frac{H}{t}-(\beta-1)(q+f_t)$ and $(\Delta f)_t=(-F^2(\nabla f)+q+f_t)_t=\frac{H}{t^2}-\frac{H_t}{t}-(\beta-1)(q_t+f_{tt})$, by Lemma \ref{lem-com-tNtD}. Thus, it is deduced from \eqref{Delta h-5} by considering the curvature condition that 
	\begin{eqnarray}\label{Delta h-6}
		\begin{split}
			\Delta^{\nabla u}H\geq&\frac{2t}{N}(F^2(\nabla f)-f_t-q)^2-2df(\nabla^{\nabla u} H)-(\beta-1)t(F^2(\nabla f))_t\\
			&-2(\beta-1)tdf(\nabla^{\nabla u} q)-2KtF^2(\nabla f)-\beta\frac{H}{t}+\beta H_t\\
			&+\beta(\beta-1)t(f_{tt}+q_t)-\beta t \Delta^{\nabla u}q.
		\end{split}
	\end{eqnarray}
	Noticing $H_t=\frac{H}{t}+t[(F^2(\nabla f))_t-\beta f_{tt}-\beta q_t]$
	, the lemma follows from \eqref{Delta h-6}. 
\end{proof}
A weak version of Lemma \ref{lem-H} is the following corollary.
\begin{cor}
	Let $(M,F,\mu)$ be a Finsler metric measure space with weighted Ricci curvature $Ric^N$ having lower bound $-K$, where $K\geq0$ is a nonnegative constant. For any $\beta\geq 1$, $H=t(F^2(\nabla f)-\beta f_t-\beta q)$ satisfies
	\begin{eqnarray}\label{equ-lem-H-int}
		\begin{split}
			&\int_0^T\int_M[d\varphi(\nabla^{\nabla u}H)-H\varphi_t]d\mu dt-\beta\int_0^T\int_Mtd\varphi(\nabla^{\nabla u}q)d\mu dt\\
			\leq&2\int_0^T\int_M\varphi[df(\nabla^{\nabla u} H)+\frac{H}{2t}+KtF^2(\nabla f)-\frac{t}{N}(F^2(\nabla f)-f_t-q)^2\\
			&\quad\quad\quad\quad\quad+(\beta-1)tdf(\nabla^{\nabla u} q)]d\mu dt,
		\end{split}
	\end{eqnarray}
	for any positive test function $\varphi\in C_0^2([0,T])\cap C^2_0(M)$, provided $f$ satisfying \eqref{equ-f}.
\end{cor}

\begin{rem}
	\begin{itemize}
		\item[i)] On the whole manifold $M$, we only could prove \eqref{Delta h-4} in the distributional sense. But we just express all the equalities in the case that $u$ and $q$ have the enough regularities. Indeed, applying the same method in the proof of integrated Bochner-Weitzenb\"ock formula (Theorem 3.6) in \cite{Ohta2014}, we could choose smooth function sequences $u_k$ and $q_k$, such that $u_k\in C^{\infty}(M)$ with $u_k\rightarrow u$ locally in the $H^2$-norm as $k\rightarrow\infty$, and simultaneously, $q_k\in C^{\infty}(M)$ with $q_k\rightarrow q$ locally in the $H^1$-norm as $k\rightarrow\infty$. We can assume $u_k$ and $q_k$ are both smooth in the variable $t$ by using Lemma \ref{lem-Ju} on $u_k$ and $q_k$ to replace them by $u_{k,\epsilon}$ and $q_{k,\epsilon}$ if necessary. We still denote them by $u_k$ and $q_k$. Furthermore, we denote $f_k=\log u_k$ and $H_k=t(F^2(\nabla^{\nabla u_k} f_k)-\beta \partial_tf_k-\beta q_k)$ to replace $f$ and $H$, respectively in the proof. Let $k$ tend to the infinity can strictly prove \eqref{Delta h-4}.
		\item[ii)] The detail to prove inequality \eqref{Delta h-6} contains an approximate sequence of positive test functions and a simple measure analysis on $M\setminus M_u$. Rigorous arguments can be made using the methods in the proof of Lemma 3.2 in \cite{Xia2023}.
	\end{itemize}	
\end{rem}

The main task in this subsection is to give a global gradient estimate of the positive solution to \eqref{equ-FS1}. We now begin with the definition below, which may require less regularity about $q$.
\begin{defn}
	We call a function $h\in H^1(M)$ on a Finsler metric measure space $(M,F,\mu)$ satisfies $\Delta^V h\leq \theta(x)$ for any local vector field $V$ in the distributional sense, if for any positive function $\varphi\in C_0^2(M)$, 
	\begin{eqnarray}
		-\int_Md\varphi(\nabla^V h)d\mu\leq \int_M\varphi\theta d\mu.
	\end{eqnarray}
\end{defn}

The Li-Yau gradient estimates of the positive solution to the Schr\"odinger equation \eqref{equ-FS1} can be obtained from Lemma \ref{lem-H} directly as follows.

\begin{thm}\label{thm-cpt-1}
	Let $(M,F,\mu)$ be a compact Finsler metric measure space with or without boundary, whose weighted Ricci curvature $Ric_N$ bounded from below by $-K$, with a nonnegative constant $K\geq0$. Suppose $q(x,t)\in C([0,T],H^1(M))$ is a function on $M\times [0,T]$. Moreover it satisfies $F_V(\nabla^V q(x))\leq\gamma$ and $\Delta^V q(x)\leq \theta$ in the distributional sense, for any point $x$ on $M$ and any local vector field $V$ around $x$, with two constants $\theta$ and $\gamma$. Assume $u(x,t)$ is a positive solution on $M\times(0,T]$ to the Finslerian Schr\"odinger equation \eqref{equ-FS1}. Moreover, $u$ satisfies the Neumann boundary condition if $\partial M\neq\emptyset$. It satisfies that
	\begin{eqnarray}\label{equ-thm1cpt}
		\begin{split}
			\sup_{M\times [0,T]}&\left\{\frac{F^2(\nabla u)}{u^2}-\frac{\beta u_t}{u}-\beta q,\quad\frac{F^2(-\nabla u)}{u^2}-\frac{\beta u_t}{u}-\beta q\right\}\\
			&\leq \frac{N\beta^2}{2t}+\left[\frac34\left(\frac{(\beta-1)^2\beta^8\gamma^4}{4\epsilon}\right)^{\frac13}+\frac{N^2\beta^4K^2}{4(1-\epsilon)(\beta-1)^2}+\frac{N}{2}\beta^3\theta\right]^{\frac12},
		\end{split}
	\end{eqnarray}
	for any $\beta>1$ and $\epsilon\in(0,1)$.
\end{thm}

\begin{proof}
	Suppose $u$ is a positive solution to the Finslerian Schr\"odinger equation, and $f=\log u$. 
	It satisfies from \eqref{equ-lem-H-int} that
	\begin{eqnarray}\label{equ-lem-H-int-2}
		\begin{split}
			&\int_0^T\int_M[d\varphi(\nabla^{\nabla u}H)-H\varphi_t-2\varphi df(\nabla^{\nabla u} H)]d\mu dt\\
			\leq&2\int_0^T\int_Mt\varphi\left[KF^2(\nabla f)-\frac{(F^2(\nabla f)-f_t-q)^2}{N}+(\beta-1)df(\nabla^{\nabla u} q)\right]d\mu dt\\
			&+\int_0^T\int_M\frac{\varphi H}{t}d\mu dt+\beta\int_0^T\int_Mtd\varphi(\nabla^{\nabla u}q)d\mu dt,
		\end{split}
	\end{eqnarray}
	for any positive test function $\varphi\in C_0^2([0,T])\cap C^2_0(M)$. By setting $y=F^2(\nabla f)$ and $z=f_t+q$, with $\epsilon=2(\frac{\beta-1}{\beta}-\frac{(\beta-1)^2}{\beta^2})=\frac{2(\beta-1)}{\beta^2}$, term $(F^2(\nabla f)-f_t-q)^2=(y-z)^2$ is equal to
	\begin{eqnarray}\label{H-inequ-2-2}
		\begin{split}
			(y-z)^2&=\frac{1}{\beta^2}[(\beta-1)y+(y-\beta z)]^2\\
			&=\frac{(\beta-1)^2}{\beta^2}y^2+\epsilon y(y-\beta z)+\frac{1}{\beta^2}(y-\beta z)^2.
		\end{split}
	\end{eqnarray} 
	Thus, the first term on the RHS of \eqref{equ-lem-H-int-2} becomes 
	\begin{eqnarray*}
		\begin{split}
			&\int_0^T\int_Mt\varphi\left[KF^2(\nabla f)-\frac{(F^2(\nabla f)-f_t-q)^2}{N}+(\beta-1)df(\nabla^{\nabla u} q)\right]d\mu dt\\
			\leq&\int_0^T\int_Mt\varphi\left[Ky-\frac{(y-\beta z)^2}{N}+(\beta-1)\gamma y^{\frac12}\right]d\mu dt\\
			=&\int_0^T\int_Mt\varphi\left[Ky+(\beta-1)\gamma y^{\frac12}-\frac1N\left(\frac{1}{\beta^2}(y-\beta z)^2+\epsilon y(y-\beta z)+\frac{(\beta-1)^2}{\beta^2}y^2\right)\right]d\mu dt\\
			\leq&\int_0^T\int_Mt\varphi\left[(\beta-1)\gamma y^{\frac12}-\frac1N\left(\frac{1}{\beta^2}(y-\beta z)^2-\frac{N^2K^2\beta^2}{4(1-\epsilon)(\beta-1)^2}+\epsilon\frac{(\beta-1)^2}{\beta^2}y^2\right)\right]d\mu dt\\
			&-\frac{\epsilon}{N}\int_0^T\int_Mt\varphi y(y-\beta z) d\mu dt\\
			\leq&\int_0^T\int_Mt\varphi\left[\frac34\left(\frac{N\beta^2(\beta-1)^2\gamma^4}{4\epsilon}\right)^{\frac13}-\frac1N\left(\frac{1}{\beta^2}(y-\beta z)^2-\frac{N^2K^2\beta^2}{4(1-\epsilon)(\beta-1)^2}\right)\right]d\mu dt\\
			&-\frac{\epsilon}{N}\int_0^T\int_Mt\varphi y(y-\beta z) d\mu dt,
		\end{split}
	\end{eqnarray*}	
	where the second inequality follows from the Young's inequality, and the third inequality follows from
	\begin{eqnarray}\label{H-inequ-4}
		\epsilon\frac{(\beta-1)^2}{\beta^2}y^2-N(\beta-1)\gamma y^{\frac12}
		\geq-\frac34\left(\frac{N^4\beta^2(\beta-1)^2\gamma^4}{4\epsilon}\right)^{\frac13}.
	\end{eqnarray}  
	Now \eqref{equ-lem-H-int-2} becomes 
	\begin{eqnarray*}
		\begin{split}
			&\int_0^T\int_M[d\varphi(\nabla^{\nabla u}H)-H\varphi_t-2\varphi df(\nabla^{\nabla u} H)]d\mu dt\\
			\leq&2\int_0^T\int_Mt\varphi\left[\frac34\left(\frac{N\beta^2(\beta-1)^2\gamma^4}{4\epsilon}\right)^{\frac13}-\frac1N\left(\frac{1}{\beta^2}(y-\beta z)^2-\frac{N^2K^2\beta^2}{4(1-\epsilon)(\beta-1)^2}\right)\right]d\mu dt\\
			&-\frac{2\epsilon}{N}\int_0^T\int_M\varphi yH d\mu dt+\int_0^T\int_M\frac{\varphi H}{t}d\mu dt+\beta\int_0^T\int_Mtd\varphi(\nabla^{\nabla u}q)d\mu dt\\
			=&2\int_0^T\int_M\frac{\varphi}{t}\left\{\frac{H}{2}+\frac{t^2}{N}\left[\frac34\left(\frac{N^4\beta^2(\beta-1)^2\gamma^4}{4\epsilon}\right)^{\frac13}+\frac{\beta^2N^2K^2}{4(1-\epsilon)(\beta-1)^2}+\frac{\beta N\theta}{2}\right]\right.\\
			&\quad\quad\quad\quad\quad\quad\left.-\frac{H^2}{N\beta^2}\right\}d\mu dt-\frac{2\epsilon}{N}\int_0^T\int_M\varphi HF^2(\nabla f) d\mu dt,
		\end{split}
	\end{eqnarray*}
	so that 
	\begin{eqnarray}\label{inequ-Xi}
		\int_0^T\int_M[d\varphi(\nabla^{\nabla u}H)-H\varphi_t]d\mu dt\leq 2\int_0^T\int_M\varphi\Xi(x,t) d\mu dt,
	\end{eqnarray}
	with 
	\begin{eqnarray*}
		\begin{split}
			\Xi=&\frac{1}{t}\left\{\frac{H}{2}+\frac{t^2}{N}\left[\frac34\left(\frac{N^4\beta^2(\beta-1)^2\gamma^4}{4\epsilon}\right)^{\frac13}+\frac{\beta^2N^2K^2}{4(1-\epsilon)(\beta-1)^2}+\frac{\beta N\theta}{2}\right]-\frac{H^2}{N\beta^2}\right\}\\
			&+df(\nabla^{\nabla u} H)-\frac{\epsilon}{N} HF^2(\nabla f).
		\end{split}	
	\end{eqnarray*}
	Suppose that $H$ achieves its maximum on $M\times [0,T]$ at $(x_0,t_0)$. Therefore, $\nabla^{\nabla u} H=0$ at $(x_0,t_0)$. If $H$ is nonpositive at $(x_0,t_0)$, the result \eqref{equ-thm1cpt} is trivial. So we suppose $H>0$ at $(x_0,t_0)$ without losing the generality. We claim that $\Xi$ is nonnegative at $(x_0,t_0)$. If not, there is a neighborhood $U$ of $(x_0,t_0)$ in $M\times(0,T]$ such that $\Xi<0$. By choosing $\varphi$ with compact support $V$ in $U$, we know from \eqref{inequ-Xi} that $H$ is a local weak supersolution of $(\partial_t-\Delta^{\nabla u})H\leq 0$. Thus, $H$ must attain the maximum on some boundary point of $V$, which is a contradiction since $(x_0,t_0)$ is an inner point of $V$.
	%
	The nonnegative of $\Xi$ at $(x_0,t_0)$ implies that 
	\begin{eqnarray*}\label{int-cpt-3}
		\frac{2}{N\beta^2}H^2-H-\frac{2t_0^2}{N}\left[\frac34\left(\frac{N^4\beta^2(\beta-1)^2\gamma^4}{4\epsilon}\right)^{\frac13}+\frac{\beta^2N^2K^2}{4(1-\epsilon)(\beta-1)^2}+\frac{\beta N\theta}{2}\right]\leq 0,
	\end{eqnarray*}
	whose LHS is a function about $H$. The large root of the function gives an upper bound of $H$ by  
	\begin{eqnarray*}
		\begin{split}
			H\leq&\frac{N\beta^2}{4}+\frac{1}{2}\left\{\frac{N^2\beta^4}{4}-4t_0^2\left[\frac34\left(\frac{N^4\beta^8(\beta-1)^2\gamma^4}{4\epsilon}\right)^{\frac13}+\frac{\beta^4N^2K^2}{4(1-\epsilon)(\beta-1)^2}+\frac{\beta^3N\theta}{2}\right]\right\}^{\frac12}\\
			\leq&\frac{N\beta^2}{2}+t_0\left[\frac34\left(\frac{N^4\beta^8(\beta-1)^2\gamma^4}{4\epsilon}\right)^{\frac13}+\frac{\beta^4N^2K^2}{4(1-\epsilon)(\beta-1)^2}+\frac{\beta^3N\theta}{2}\right]^{\frac12},
		\end{split}
	\end{eqnarray*}
	for any $\epsilon\in(0,1)$, which is equal to the first part of \eqref{equ-thm1cpt}. 
	
	For the reverse metric $\overleftarrow{F}$, it satisfies that $\overleftarrow{Ric}^N(y)=Ric^N(-y)$. Consequently, the weighted Ricci curvature $\overleftarrow{Ric}^N(y)$ of $(M,\overleftarrow{F},\mu)$ has the same lower bound $-K$ as the $Ric^N$ of $(M,F,\mu)$. The above arguments are still valid for $\overleftarrow{F}$.   
	
	The estimation of the Neumann positive solution is the same. Combining all these estimates can finish the proof.
\end{proof}

\subsection{Weighted Ricci curvature $Ric^{\infty}$ bounded from below}\label{Ch-GGE-2}

With the same method, we will show another global gradient estimates of the positive solution to the Finslerian Schr\"odinger equation \eqref{equ-FS1} on manifolds satisfying the different curvature condition, that is, the weighted Ricci curvature $Ric^{\infty}$ bounded from below by a function $-K$.

Let $u$ be a positive solution to the Finslerian Schr\"odinger equation \eqref{equ-FS1} with upper bound $L$.
We denote $f=\log\frac{u}{L}$ and $W=\frac{F^2(\nabla f)}{(1-f)^2}$. As Lemma \ref{lem-f}, we have 
\begin{lem}
	Let $u$ be a positive solution on $M\times [0,T]$ to \eqref{equ-FS1}. $f(x,t)=\log\frac{u}{L}$ satisfies that 
	\begin{eqnarray}
		f_t=\Delta f+F^2{\nabla f}-q,
	\end{eqnarray}
	on $M_u\times[0,T]$. Moreover, 
	\begin{eqnarray}\label{equ-f2}
		\int_0^T\int_M[d\varphi(\nabla f)-\varphi_tf]d\mu dt=\int_0^T\int_M\varphi(F^2(\nabla f)-q)d\mu dt,
	\end{eqnarray}
	for any test function $\varphi\in C^2(M\times [0,T])$.
\end{lem}

The next lemma is also established on $M_u\times [0,T]$.

\begin{lem}\label{lem-Wineq}
	Let $(M,F,\mu)$ be a Finsler metric measure space with weighted Ricci curvature $Ric^{\infty}$ having lower bound $-K$, where $K\geq0$ is a nonnegative constant. $W=\frac{F^2(\nabla f)}{(1-f)^2}$ satisfies
	\begin{eqnarray}\label{equ-lem-W-int}
		\begin{split}
			\frac12(\partial_t-\Delta^{\nabla u})W\leq K\frac{F^2(\nabla f)}{(1-f)^2}-\frac{dq(\nabla f)}{(1-f)^2}-\frac{fdW(\nabla f)}{1-f}+\frac{qW}{1-f}-(1-f)W^2,
		\end{split}
	\end{eqnarray}
	on $M_u\times[0,T]$, provided $f=\log\frac{u}{L}$.
\end{lem}
\begin{proof}
	From the definition of $W$, one can obtain the following computations by noticing the vanishing of the Cratan tensor in the direction $\nabla u$. That is, 
	\begin{eqnarray}\label{6.3}
		\begin{split}
			W_t=&2\frac{g_{\nabla u}(\nabla f,\nabla f_t)}{(1-f)^2}+2\frac{F^2(\nabla f)f_t}{(1-f)^3}\\
			=&\frac{2}{(1-f)^2}[df(\nabla^{\nabla u} \Delta f)+df(\nabla^{\nabla u}F^2(\nabla f))-dq(\nabla f)]\\
			&+2\frac{F^2(\nabla f)}{(1-f)^3}(\Delta f+F^2(\nabla f)-q),
		\end{split}
	\end{eqnarray}
	and
	\begin{eqnarray}\label{6.4}
		\nabla^{\nabla u}W=\frac{\nabla^{\nabla u}(F^2(\nabla f))}{(1-f)^2}+2\frac{F^2(\nabla f)\nabla f}{(1-f)^3}.
	\end{eqnarray}
	Taking another derivative of \eqref{6.4} and taking the trace with respect to $g_{\nabla u}$ yield 
	\begin{eqnarray}\label{6.5}
		\begin{split}
			\mathrm{tr}_{\nabla u}\nabla^{\nabla u2}W=&2\frac{df(\mathrm{tr}_{\nabla u}\nabla^{\nabla u2}\nabla f)+\|\nabla^2f\|_{HS(\nabla u)}^2}{(1-f)^2}+4\frac{df(\nabla^{\nabla u}F^2(\nabla f))}{(1-f)^3}\\
			&+2\frac{F^2(\nabla f)\mathrm{tr}_{\nabla u}\nabla^2f}{(1-f)^2}+6\frac{F^4(\nabla f)}{(1-f)^4}\\
			=&2\frac{df(\nabla^{\nabla u}(\mathrm{tr}_{\nabla u}\nabla^2f))+Ric(\nabla f)+\|\nabla^2f\|_{HS(\nabla u)}^2}{(1-f)^2}\\
			&+4\frac{df(\nabla^{\nabla u}F^2(\nabla f))}{(1-f)^3}+2\frac{F^2(\nabla f)\mathrm{tr}_{\nabla u}\nabla^2f}{(1-f)^2}+6\frac{F^4(\nabla f)}{(1-f)^4}.
		\end{split}
	\end{eqnarray}
	On the other hand,
	\begin{eqnarray}\label{6.6}
		\begin{split}
			\Delta^{\nabla u}W=&\mathrm{tr}_{\nabla u}\nabla^{\nabla u2}W-d\tau(\nabla^{\nabla u} W)\\
			=&\mathrm{tr}_{\nabla u}\nabla^{\nabla u2}W-\frac{d\tau(\nabla^{\nabla u}F^2(\nabla f))}{(1-f)^2}-2\frac{F^2(\nabla f)d\tau(\nabla f)}{(1-f)^3}.
		\end{split}
	\end{eqnarray}
	Combining \eqref{6.3}-\eqref{6.6}, one may find that
	\begin{eqnarray}\label{6.7}
		\begin{split}
			(\Delta^{\nabla u}&-\partial_t)W=\frac{2}{(1-f)^2}[df(\nabla^{\nabla u}(\mathrm{tr}_{\nabla u}\nabla^2f))+Ric(\nabla f)+\|\nabla^2f\|_{HS(\nabla u)}^2\\
			&-d\tau\otimes df(\nabla^2f)-df(\nabla^{\nabla u}\Delta f)-df(\nabla^{\nabla u}F^2(\nabla f))+dq(\nabla f)]\\
			&+\frac{2}{(1-f)^3}[2df(\nabla^{\nabla u}F^2(\nabla f))-F^4(\nabla f)+qF^2(\nabla f)]+6\frac{F^4(\nabla f)}{(1-f)^4}
		\end{split}
	\end{eqnarray}
	holds on $M_u\times[0,T]$.
	
	Since
	\begin{eqnarray}
		\nabla^{\nabla u}(S(\nabla f))=\dot{S}(\nabla f)+d\tau\otimes df(\nabla^2f),
	\end{eqnarray}
	\eqref{6.7} becomes 
	\begin{eqnarray}\label{6.8}
		\begin{split}
			(\Delta^{\nabla u}-\partial_t)W=&\frac{2}{(1-f)^2}[Ric^{\infty}(\nabla f)+\|\nabla^2f\|_{HS(\nabla u)}^2-2(df)^2(\nabla^2f)+dq(\nabla f)]\\
			&+\frac{2}{(1-f)^3}[4(df)^2(\nabla^2f)-F^4(\nabla f)+qF^2(\nabla f)]+6\frac{F^4(\nabla f)}{(1-f)^4}.
		\end{split}
	\end{eqnarray}
	
	Noticing that
	\begin{eqnarray}\label{6.9}
		\begin{split}
			&\frac{\|\nabla^2f\|_{HS(\nabla u)}^2}{(1-f)^2}+2\frac{(df)^2(\nabla^2f)}{(1-f)^3}+2\frac{F^4(\nabla f)}{(1-f)^4}\\
			=&\|\frac{\nabla^2 f}{1-f}+\frac{\nabla f\otimes \nabla f}{(1-f)^2}\|^2_{HS(\nabla u)}\geq 0.
		\end{split}
	\end{eqnarray}	
	It follows from
	\begin{eqnarray*}
		dW(\nabla f)=2\frac{(df)^2(\nabla^2f)}{(1-f)^2}+2\frac{F^4(\nabla f)}{(1-f)^3}
	\end{eqnarray*}	
	that
	\begin{align}\label{6.10}
		2\frac{dW(\nabla f)}{1-f}-4\frac{F^4(\nabla f)}{(1-f)^4}&=\frac{4g_{\nabla u}(\nabla^2f,\nabla f\otimes \nabla f)}{(1-f)^3},\\
		-2dW(\nabla f)+4\frac{F^4(\nabla f)}{(1-f)^3}&=-4\frac{g_{\nabla u}(\nabla^2f,\nabla f\otimes \nabla f)}{(1-f)^2}.
	\end{align}
	Plugging \eqref{6.9}-\eqref{6.10} into \eqref{6.8}, one can deduce that
	\begin{eqnarray}
		\begin{split}
			(\Delta^{\nabla u}-\partial_t)W\geq&2\frac{Ric^{\infty}(\nabla f)}{(1-f)^2}-4\frac{(df)^2(\nabla^2f)}{(1-f)^2}\\
			&+2\frac{dq(\nabla f)}{(1-f)^2}+4\frac{g_{\nabla u}(\nabla^2f,\nabla f\otimes \nabla f)}{(1-f)^3}\\
			&-2\frac{F^4(\nabla f)}{(1-f)^3}+2\frac{qF^2(\nabla f)}{(1-f)^3}+4\frac{F^4(\nabla u)}{(1-f)^4}\\
			=&2\frac{Ric^{\infty}(\nabla f)}{(1-f)^2}+2\frac{dq(\nabla f)}{(1-f)^2}\\
			&+\frac{2f}{1-f}dW(\nabla f)+2\frac{qW}{1-f}+2(1-f)W^2,
		\end{split}
	\end{eqnarray}
	On $M_u\times [0,T]$. Then the lemma can be obtained by considering the curvature condition $Ric^{\infty}\geq-K$.	
\end{proof}	

In the remaining part of this section, we are going to prove Theorem \ref{thm-cpt-2}. 

\begin{thm}\label{thm-cpt-2}
	Let $(M,F,\mu)$ be a compact Finsler metric measure space with or without boundary, whose weighted Ricci curvature $Ric^{\infty}$ is bounded from below by $-K$, with a nonnegative function $K\geq0$. Suppose $q(x,t)\in C((-\infty,\infty),H^1(M))$ is a function on $M\times [0,T]$. Moreover it satisfies $F_V(\nabla^V q(x))\leq\gamma$, for any point $x$ on $M$ and any local vector field $V$ around $x$, with nonnegative function $\gamma\geq0$. Assume $u(x,t)\leq L$ is a positive solution on $M\times[0,T]$ to \eqref{equ-FS1}. Moreover, $u$ satisfies the Neumann boundary condition if $\partial M\neq\emptyset$. It satisfies that
	\begin{eqnarray*}
		\sup_{M\times [0,T]}\left\{\frac{F(\nabla u)}{u},\quad\frac{F(-\nabla u)}{u}\right\}\leq \sqrt{2}(\|K\|_{L^{\infty}}^{\frac12}+\|q^-\|_{L^{\infty}}^{\frac12}+\|\gamma\|_{L^{\infty}}^{\frac13})(1+\log\frac{L}{u}).
	\end{eqnarray*}
\end{thm}

\begin{proof}
	From \eqref{equ-lem-W-int}, we deduce that                
	\begin{eqnarray}\label{6.11}
		\begin{split}
			&\frac12\int_0^T\int_M\left\{d\varphi(\nabla^{\nabla u}W)-W\varphi_t+2\varphi[\frac{ f}{1-f}dW(\nabla f)+(1-f)W^2]\right\}d\mu dt\\
			\leq&\int_0^T\int_M\varphi\left[K\frac{F^2(\nabla f)}{(1-f)^2}-\frac{dq(\nabla f)}{(1-f)^2}+\frac{qW}{1-f}\right]d\mu dt,
		\end{split}
	\end{eqnarray}
	for any positive test function $\varphi\in C_0^2([0,T])\cap C^2_0(M)$, provided $f$ satisfying \eqref{equ-f}.
	
	According to the Young's inequality, we have the following inequalities.
	\begin{eqnarray}\label{6.12}
		\begin{split}
			\int_0^T\int_M\varphi K\frac{F^2(\nabla f)}{(1-f)^2}&d\mu dt=\int_0^T\int_M\varphi KWd\mu dt\\
			&\leq \frac14\int_0^T\int_M\varphi(1-f)W^2d\mu dt+\int_0^T\int_M\frac{\varphi K^2}{1-f}d\mu dt,
		\end{split}
	\end{eqnarray}
	and
	\begin{eqnarray}
		\begin{split}
			-\int_0^T\int_M\varphi\frac{qW}{1-f}d\mu dt&\leq\int_0^T\int_M\varphi\frac{Wq^-}{1-f}d\mu dt\\
			&\leq \frac14\int_0^T\int_M\varphi(1-f)W^2d\mu dt+\int_0^T\int_M\frac{\varphi(q^-)^2}{(1-f)^3}d\mu dt,
		\end{split}	
	\end{eqnarray}
	where $q^-=\max\{-q,0\}$. Moreover, by the same reason, we have
	\begin{eqnarray}\label{6.13}
		\begin{split}
			-\int_0^T\int_M\varphi&\frac{dq(\nabla f)}{(1-f)^2}d\mu dt\leq \int_0^T\int_M\varphi\frac{F^*_{\nabla u}(dq)}{1-f}W^{\frac12}d\mu dt\\
			&\leq\frac14\int_0^T\int_M\varphi(1-f)W^2d\mu dt+\frac{3}{4}\int_0^T\int_M\frac{\varphi F^{\frac43}_{\nabla u}(\nabla q)}{(1-f)^{\frac53}}d\mu dt.			
		\end{split}
	\end{eqnarray}
	Plugging inequalities \eqref{6.12}-\eqref{6.13} into \eqref{6.11} shows that
	\begin{eqnarray}\label{6.14}
		\begin{split}
			\frac12\int_0^T\int_M[d\varphi(\nabla^{\nabla u}W)&-W\varphi_t]d\mu dt\leq\int_0^T\int_M\varphi\left[\frac{K^2}{1-f}+\frac{(q^-)^2}{(1-f)^3}\right.\\
			&\left.+\frac{3}{4}\frac{F^{\frac43}_{\nabla u}(\nabla q)}{(1-f)^{\frac53}}-\frac{f}{1-f}dW(\nabla f)-\frac{(1-f)}{4}W^2\right]d\mu dt,
		\end{split}
	\end{eqnarray}
	By the same argument in the proof of Theorem \ref{thm-cpt-1}, we denote the RHS of \eqref{6.14} by $\int_0^T\int_M\varphi\Xi d\mu dt$ with $\Xi=\frac{K^2}{1-f}+\frac{(q^-)^2}{(1-f)^3}+\frac{3}{4}\frac{F^{\frac43}_{\nabla u}(\nabla q)}{(1-f)^{\frac53}}-\frac{f}{1-f}dW(\nabla f)-\frac{(1-f)}{4}W^2$. We claim that $\Xi\geq0$ at the maximum point $(x_1,t_1)$ of $W$, otherwise $W$ should be a local weak supersolution to $(\partial_t-\Delta^{\nabla u})W\leq 0$, deducing a contradiction since the $(x_1,t_1)$ is a inner maximum point of $W$.
	The inequality $\Xi\geq 0$ is equal to the following at $(x_1,t_1)$,
	\begin{eqnarray}
		\frac{F^2(\nabla f)}{(1-f)^2}\leq\frac{2K}{1-f}+\frac{2q^-}{(1-f)^2}+\sqrt{3}\frac{F^{\frac23}_{\nabla u}(\nabla q)}{(1-f)^{\frac43}}.
	\end{eqnarray}
	Since $f\leq 0$, we finally derive that
	\begin{eqnarray}
		\frac{F(\nabla f)}{1-f}\leq\sqrt{2}(\|K\|_{L^{\infty}}^{\frac12}+\|q^-\|_{L^{\infty}}^{\frac12}+F_{\nabla u}^{\frac13}(\nabla q)).
	\end{eqnarray}
	The above arguments are still valid for $\overleftarrow{F}$, and the estimation of the Neumann positive solution is the same. Thus we finish the proof.  	
\end{proof}

\begin{rem}
	The condition of $F_V(\nabla^V q(x))\leq \gamma$ in both Theorem \ref{thm-cpt-1} and Theorem \ref{thm-cpt-2} can be repalaced by $F^*(dq(x))=F(\nabla q)\leq \gamma$. It is because $F_V(\nabla^V q(x))=F_{V}^*(dq(x))$. Since the manifold $M$ is compact, it admits finite uniform convexity and finite uniform smoothness, hence has finite misalignment. Based on this fact, we will take $F(\nabla q)\leq \gamma$ to describe the bound of $\nabla^{\nabla u} q$ in the rest of this article, since we always assume there is a finite misalignment within the area we look into.
\end{rem}

\section{Local gradient estimates on forward complete Finsler manifolds with two kinds of mixed weighted Ricci curvature bounded from below}\label{Ch-noncpt}

In this section, we intend to discuss gradient estimation of the Finslerian Schr\"odinger equation further on forward complete Finsler metric measure spaces. The gradient estimate on a noncompact Finsler manifold is a hard topic, which was first asked by Ohta in \cite{Ohta2014}. Later, Q. Xia answered the question partially in \cite{Xia2023} by using the Moser's iteration method, which was first proposed by C. Xia \cite{XiaC2014CVPDE}. Actually, Q. Xia got the gradient estimates of the positive solution to the Finslerian heat equation on forward complete Finsler metric measure spaces in the weighted Ricci curvature condition. However, the method in \cite{Xia2023} is not universal, even for the positive solution to the Finslerian heat equation, when the estimation term with $\beta=1$ (cf. Remark 4.1 in \cite{Xia2023}).

The main obstacle to reach the gradient estimates on forward complete noncompact Finsler manifolds is the absence of the Laplacian comparison theorem of the distance function. To overcome it, we introduce a new curvature tensor called the mixed weighted Ricci curvature, as well as some non-Riemannian tensors (cf. Subsection \ref{subsec-2.1}) to get the Laplacian comparison theorem (cf. Section \ref{sec-Lapcomp}). This method can be applied to lots of PDEs on forward complete Finsler metric measure spaces.

As we illustrate in the last section, the derivatives about $x$ and $t$ in this section are also valid for the local positive solutions to the Finslerian Schr\"odinger equation \eqref{equ-FS1}, by uitilizing the continuous version and the smooth approximation $u_{\epsilon}$ in $t$.

With the assistance of the Laplacian comparison theorem \ref{thm-LapComp-1} and the computation in the last section, we are ready to prove Theorem \ref{thm-GEcomplete-1} now.

\begin{proof}[Proof of Theorem \ref{thm-GEcomplete-1}]
	As before, we define $H=t(F^2(\nabla f)-\beta f_t-\beta q)$ with $f=\log u$. Let $\tilde \phi(r)$ is a smooth function satisfies
	\begin{eqnarray}
		\tilde\phi(r)=\begin{cases}
			1\quad r\in[0,1]\\
			0\quad r\in[2,\infty),
		\end{cases}
	\end{eqnarray}
	$-C_1\leq\frac{\tilde\phi'(r)}{\sqrt{\tilde \phi(r)}}\leq 0$, and $\tilde\phi''\geq -C_2$, where $C_1,C_2$ are two positive constants. We set the cut-off function as $\phi(x)=\tilde\phi(\frac{r(x)}{R})$, and consider the function $\phi H$ on $B(p,2R)\times[0,\infty)$. 
	
	Noticing that $\Delta^{\nabla u}(\phi H)=(\Delta^{\nabla u}\phi)H+2g_{\nabla u}(\nabla^{\nabla u} \phi,\nabla^{\nabla u} H)+\phi(\Delta^{\nabla u}H)$, and following the same test function method in Lemma \ref{lem-H}, we could get the following inequality
	\begin{eqnarray}\label{DeltaphiH-1}
		\begin{split}
			\Delta^{\nabla u}(\phi H)\geq&-H\Delta^{\nabla u}\phi+2g_{\nabla u}(\nabla^{\nabla u}\log\phi,\nabla^{\nabla u}(\phi H))-2HF_{\nabla u}^2(\nabla^{\nabla u}\phi)\phi^{-1}\\
			&+\phi\left[H_t-2df(\nabla^{\nabla u} H)+\frac{2t}{N}(F^2(\nabla f)-f_t-q)^2-2KtF^2(\nabla f)\right.\\
			&\quad\quad\,\,\left.-\frac{H}{t}-\beta t\Delta^{\nabla u}q-2(\beta-1)tdf(\nabla^{\nabla u} q)\right],
		\end{split}
	\end{eqnarray}
	which implies that
	\begin{eqnarray}\label{DeltaphiH-2}
		\begin{split}
			&\int_0^T\int_{B_{2R}}\{-d\varphi(\nabla^{\nabla u}(\phi H))+\varphi_t\phi H\} d\mu dt\\
			\geq&\int_0^T\int_{B_{2R}}\varphi\{H[\Delta^{\nabla u}\phi-\frac{2}{\phi}F_{\nabla u}^2(\nabla^{\nabla u}\phi)]-\frac{\phi H}{t}+2d\log\phi(\nabla^{\nabla u}(\phi H))\\
			&\quad\quad\quad\quad\quad\,\,-\frac{2C_1}{R}\vert H\vert F(\nabla f)\phi^{\frac12}+\frac{2\phi t}{N}(F^2(\nabla f)-f_t-q)^2-2Kt\phi F^2(\nabla f)\\
			&\quad\quad\quad\quad\quad\,\,-\beta t\phi \theta-2(\beta-1)t\phi\gamma F(\nabla f)+2df(\nabla^{\nabla u}(\phi H))\}d\mu dt,
		\end{split}
	\end{eqnarray}
	for any $C^2$ smooth positive function $\varphi$ on $M\times[0,T]$. Denoting the RHS of \eqref{DeltaphiH-2} by $\int_0^T\int_{B_{2R}}\varphi\Xi(x,t)d\mu dt$, we could deduce $\Xi\leq 0$ at $(x_2,t_2)$, where $\phi H$ attains its positive maximum in $B(p,2R)$, by the same argument as in the compact cases (cf. the proofs of Theorem \ref{thm-cpt-1} and Theorem \ref{thm-cpt-2}). Therefore, at $(x_2,t_2)$, we have $\nabla^{\nabla u}(\phi H)=0$ and 		
	\begin{eqnarray}\label{DeltaphiH-4}
		\begin{split}
			0\geq& H[tJ(\phi,R)-1]-t^2\beta\theta\\
			&+\frac{2t^2}{N}\left\{(y-z)^2-NKy-N(\beta-1)\gamma y^{\frac12}-NC_1R^{-1}y^{\frac12}(y-\beta z)\right\},
		\end{split}
	\end{eqnarray}
	where $J(\phi,R)=\Delta^{\nabla u}\phi-2F_{\nabla u}^2(\nabla^{\nabla u}\phi)\phi^{-1}$, $y=\phi F^2(\nabla f)$ and $z=\phi(f_t+q)$.

	Setting $\epsilon=\frac{2(\beta-1)}{\beta^2}$ again and employing \eqref{H-inequ-2-2} yields 
	\begin{eqnarray*}
		\begin{split}
			& (y-z)^2-NKy-N(\beta-1)\gamma y^{\frac12}-NC_1R^{-1}y^{\frac12}(y-\beta z)\\
			=&\frac{1}{\beta^2}(y-\beta z)^2+(\epsilon y-NC_1R^{-1}y^{\frac12})(y-\beta z)\\
			&+\frac{(\beta-1)^2}{\beta^2}y^2-NKy-N(\beta-1)\gamma y^{\frac12},
		\end{split}
	\end{eqnarray*}
	in which 
	$$\epsilon y-NC_1R^{-1}y^{\frac12}=2(\beta-1)\beta^{-1}y-NC_1R^{-1}y^{\frac12}\geq-\frac{N^2}{8}C_1^2\frac{\beta^2}{\beta-1}R^{-2},$$
	so that  
	\begin{eqnarray}\label{DeltaphiH-5}
		\begin{split}
			&\quad (y-z)^2-NKy-N(\beta-1)\gamma y^{\frac12}-NC_1R^{-1}y^{\frac12}(y-\beta z)\\
			&\geq\frac{1}{\beta^2}(y-\beta z)^2-\frac{C_1N^2\beta^2}{8(\beta-1)R^2}(y-\beta z)+\frac{(\beta-1)^2}{\beta^2}y^2-NKy-N(\beta-1)\gamma y^{\frac12}.
		\end{split}
	\end{eqnarray}
	The last three terms in \eqref{DeltaphiH-5} can be estimated as
	\begin{eqnarray}
		\begin{split}
			&\frac{(\beta-1)^2}{\beta^2}y^2-NKy-N(\beta-1)\gamma y^{\frac12}\\
			\geq& \frac{(\beta-1)^2}{\beta^2}y^2-(1-\epsilon)\frac{(\beta-1)^2}{\beta^2}y^2-\frac{N^2K^2\beta^2}{4(1-\epsilon)(\beta-1)^2}-N(\beta-1)\gamma y^{\frac12}\\
			=&\epsilon\frac{(\beta-1)^2}{\beta^2}y^2-N(\beta-1)\gamma y^{\frac12}-\frac{\beta^2N^2K^2}{4(1-\epsilon)(\beta-1)^2}\\
			\geq&-\frac34(\frac{N^4\beta^2(\beta-1)^2\gamma^4}{4\epsilon})^{\frac13}-\frac{\beta^2N^2K^2}{4(1-\epsilon)(\beta-1)^2},
		\end{split}
	\end{eqnarray}
	for any $\epsilon\in(0,1)$.
	Plugging it into \eqref{DeltaphiH-5} yields
	\begin{eqnarray}\label{DeltaphiH-6}
		\begin{split}
			&(y-z)^2-NKy-N(\beta-1)\gamma y^{\frac12}-NC_1R^{-1}y^{\frac12}(y-\beta z)\\
			\geq&\frac{1}{\beta^2}(y-\beta z)^2-\frac{C_1N^2\beta^2}{8(\beta-1)R^2}(y-\beta z)\\
			&-\frac34(\frac{N^4\beta^2(\beta-1)^2\gamma^4}{4\epsilon})^{\frac13}-\frac{\beta^2N^2K^2}{4(1-\epsilon)(\beta-1)^2}.
		\end{split}
	\end{eqnarray}
	It could be deduced from the assumptions of the derivative estimates of $\tilde\phi$ that
	\begin{eqnarray*}
		\begin{split}
			&\quad\Delta^{\nabla u}\phi=\frac{\tilde\phi'\Delta^{\nabla u} r}{R}+\frac{\tilde\phi''\vert \nabla^{\nabla u} r\vert ^2}{R^2}\\
			&\geq-\frac{C_1}{R}\left[C(N,A)\sqrt{\frac{K(2R)}{C(N,A)}}\coth\left(R\sqrt{\frac{K(2R)}{C(N,A)}}\right)+C_0(K_0,A)\right]-\frac{C_2}{R^2},		
		\end{split}
	\end{eqnarray*}
	which is equal to say
	\begin{eqnarray}\label{est-J}
		\begin{split}
			J(\phi,R)\geq&-\frac{C_1}{R}\left[C_1R^{-1}\sqrt{C(N,A)K(2R)}\coth\left(R\sqrt{\frac{K(2R)}{C(N,A)}}\right)+C_0\right]\\
			&-\frac{C_2}{R^{2}}-\frac{2C_1^2}{R^{2}}\\
			\geq& -\frac{C_3}{R^{2}}[1+R(1+\sqrt{K(2R)})],
		\end{split}
	\end{eqnarray}
	where $C_3=C_3(N,A,K_0)$ is a constant depending on $K_0$, $N$ and $A$. Combining \eqref{DeltaphiH-6} and \eqref{est-J}, we can find from  \eqref{DeltaphiH-4} that
	\begin{eqnarray}
		\begin{split}
			0\geq&-\phi H[t_0C_3R^{-2}(1+R(1+\sqrt{K(2R)}))+1]-t_0^2\beta\theta\\
			&+\frac{2t_0^2}{N}[\frac{1}{\beta^2}(y-\beta z)^2-\frac{C_1N^2\beta^2}{8(\beta-1)R^2}(y-\beta z)\\
			&\quad\quad\quad-\frac34(\frac{N^4\beta^2(\beta-1)^2\gamma^4}{4\epsilon})^{\frac13}-\frac{\beta^2N^2K^2}{4(1-\epsilon)(\beta-1)^2}]\\
			&=\frac{2}{N\beta^2}(\phi H)^2-\left[\frac{C_3t_0}{R^2}(1+R(1+\sqrt{K}))+\frac{NC_1\beta^2t_0}{4(\beta-1)R^2}+1\right](\phi H)\\
			&\quad-t_0^2\left[\beta\theta+\frac32\left(\frac{N\beta^2(\beta-1)^2\gamma^4}{4\epsilon}\right)^{\frac13}+\frac{\beta^2NK^2}{4(1-\epsilon)(\beta-1)^2}\right].
		\end{split}
	\end{eqnarray}
	Equivalently
	\begin{eqnarray}
		\begin{split}
			0&\geq (\phi H)^2-\left[\frac{C_3N\beta^2t_0}{2R^2}(1+R(1+\sqrt{K}))+\frac{C_1N^2\beta^4t_0}{8(\beta-1)R^2}+\frac{N\beta^2}{2}\right](\phi H)\\
			&\quad-t_0^2\left[\frac{N\beta^3}2\theta+\frac34\left(\frac{N^4\beta^8(\beta-1)^2}{4\epsilon}\right)^{\frac13}\gamma^{\frac43}+\frac{\beta^4N^2K^2}{4(1-\epsilon)(\beta-1)^2}\right],
		\end{split}
	\end{eqnarray}
	which implies that 
	\begin{eqnarray}
		\begin{split}
			\phi H\leq &t_0\left[\frac{C_3N\beta^2t_0}{2R^2}(1+R(1+\sqrt{K}))+\frac{C_1N^2\beta^4t_0}{8(\beta-1)R^2}\right]+\frac{N\beta^2}{2}\\
			&+t_0\left[\frac{N\beta^3}2\theta+\frac34\left(\frac{N^4\beta^8(\beta-1)^2}{4\epsilon}\right)^{\frac13}\gamma^{\frac43}+\frac{\beta^4N^2K^2}{4(1-\epsilon)(\beta-1)^2}\right]^{\frac12}.
		\end{split}
	\end{eqnarray}
	the result \eqref{equ-thm1} follows by taking $\phi\equiv1$ on $B(p,R)$ and $\phi\equiv0$ out of $B(p,2R)$. 	
\end{proof}

The following theorem is a direct corollary of Theorem \ref{thm-GEcomplete-1}.

\begin{thm}\label{thm-GEcomplete-1-cor}
	Let $(M,F,\mu)$ be a forward complete Finsler metric measure space without boundary. Suppose it has finite misalignment $\alpha\leq\infty$ and the mixed weighted Ricci curvature $^mRic^N$ is bounded from below by $-K$ with a positive constant $K$. Moreover, the non-Riemannian tensors satisfy $F(U)+F^*(\mathcal{T})+F(\mathrm{div} C(V))\leq K_0$,
	for any local vector fields $V,W$ in the definitions of $U$ and $\mathcal{T}$. Suppose $u$ is a positive solution to the Schr\"odinger equation \eqref{equ-FS1}, and suppose there is a point $p\in M$ such that $F(\nabla q)\leq \gamma(r(x),t)$ and $\Delta^{V}q\leq \theta$ on $M\times (0,T]$ with a constant $\theta$, for any locally nonvanishing vector field $V$.
	\begin{itemize}
		\item[(i)] If $K=0$ and $\lim_{r\rightarrow \infty}\frac{\gamma(r,t)}{r}\leq\iota(t)$, then
		\begin{eqnarray}
			\frac{F^2(\nabla u)}{u^2}-\frac{u_t}{u}\leq\frac{N}{2t}+\sqrt{\frac{N\theta}{2}}+q+C_5\sqrt{\iota(t)}.
		\end{eqnarray}
		\item[(ii)] If $\gamma(r,t)\leq \gamma_0(t)$, then
		\begin{eqnarray}
			\frac{F^2(\nabla u)}{u^2}-\frac{\beta u_t}{u}\leq\frac{N\beta^2}{2t}+\beta q+C_6\left(\gamma_0^{\frac23}(t)+\frac{K}{\beta-1}+\theta^{\frac12}\right).
		\end{eqnarray}	
	\end{itemize}
\end{thm}

\begin{proof}
	Taking $\beta=R^{-2}\iota^{-\frac12}+1$ in \eqref{equ-thm1}, we have 
	\begin{eqnarray}
		\begin{split}
			&\quad\frac{F^2(\nabla u)}{u^2}-(R^{-2}\iota^{-\frac12}+1)\frac{u_t}{u}-(R^{-2}\iota^{-\frac12}+1)q\\
			&\leq C_3(R^{-2}\iota^{-\frac12}+1)^2R^{-2}(1+R+\frac{(R^{-2}\iota^{-\frac12}+1)^2}{R^{-2}\iota^{-\frac12}})+\frac{N}{2}(R^{-2}\iota^{-\frac12}+1)^2t^{-1}\\
			&\quad+\left[\frac{N}{2}(R^{-2}\iota^{-\frac12}+1)^3\theta\right]^{\frac12}+C_4^{\frac12}\left(\frac{(R^{-2}\iota^{-\frac12})(R^{-2}\iota^{-\frac12}+1)^4\gamma^2}{\sqrt{\epsilon}}\right)^{\frac13}.
		\end{split}
	\end{eqnarray}
	Setting $R\rightarrow\infty$ and $\epsilon=\frac12$ can provide that
	\begin{eqnarray}
		\frac{F^2(\nabla u)}{u^2}-\frac{u_t}{u}-q\leq C_3\iota^{\frac12}+\frac{N}{2t}+(\frac{N\theta}{2})^{\frac12}+(2^{\frac13}C_4)^{\frac12}\iota^{\frac12},
	\end{eqnarray}	
	which is equal to the result in (i). The result in (ii) is obtained by noticing the boundedness of $(\beta-1)^2\eta^8$, $\beta^4$ and $N\beta^3/2$ in \eqref{equ-thm1}.
\end{proof}


\begin{rem}
	We also can replace the bounds of the mixed weighted Ricci curvature and the non-Riemannian tensors in Theorem \ref{thm-GEcomplete-1} and Theorem \ref{thm-GEcomplete-1-cor} by the lower bound of weighted flag curvature according to Theorem \ref{thm-LapComp-2}.\\
\end{rem}

For completeness, we are going to end this section by the gradient estimates of the positive solution on forward complete Finsler manifold with bounded $^mRic^{\infty}$, corresponding to Theorem \ref{thm-cpt-2}. To apply the comparison theorem in the mixed weighted Ricci curvature $^mRic^k_W(V)$, We remind readers to notice Remark \ref{rem-comparison}.

\begin{proof}[Proof of Theorem \ref{thm-GEcomplete-2}]
	Let $\phi=\phi(x,t)$ be a smooth cut-off function supported in $Q_{2R,T}$, defined by $\phi(x,t)=\tilde\phi(d(p,x),t)$, where $d$ denotes the forward distance. The function $\tilde\phi$ satisfying the following properties
	\begin{itemize}
		\item[i)] $0\leq\tilde\phi\leq1$, and $\tilde\phi=1$ in $B(p,R)\times[t_0-T,t_0]$, $\tilde\phi=0$ out of $Q_{2R,T}$;
		\item[ii)] $-C\leq \tilde\phi_r/\sqrt\phi\leq 0$ and $\tilde\phi_{rr}\geq-C$;
		\item[iii)] $\vert \tilde\phi_t\vert /\sqrt\phi\leq C/T$.
	\end{itemize}
	Such function is a special one of the cut-off function explored by Souplet and Zhang \cite{SZ2006BLMS}.
	It follows from \eqref{equ-lem-W-int} that
	\begin{eqnarray}\label{7p1}
		\begin{split}		
			&\int_0^T\int_M\{-d\varphi(\nabla^{\nabla u}(\phi W))+\phi W\varphi_t\\
			&\quad\quad\quad\quad\quad-\varphi[\frac{2f}{1-f}df(\nabla^{\nabla u}(\phi W))-\frac{2}{\phi}d\phi(\nabla^{\nabla u}(\phi W))]\}d\mu dt\\
			\geq&2\int_0^T\int_M\varphi\left[(1-f)\phi W^2-\frac{f}{1-f}d\phi(\nabla f)W-\frac{F_{\nabla u}^2(\nabla^{\nabla u}\phi)}{\phi}W\right.\\
			&\left.\quad\quad\quad\quad\quad\quad+\frac12(\Delta^{\nabla u}\phi-\phi_t)W-K\phi W+\frac{q}{1-f}\phi W+\frac{\phi dq(\nabla f)}{(1-f)^2}\right]d\mu dt,
		\end{split}
	\end{eqnarray}
	for any positive test function $\varphi\in C_0^2([0,T])\cap C^2_0(M)$.
	
	Noticing $d\phi(\nabla f)\leq F_{\nabla u}(\nabla^{\nabla u}\phi)F(\nabla f)$ and $ dq(\nabla f)\leq  F_{\nabla u}(\nabla^{\nabla u}q)F(\nabla f)$, \eqref{7p1} becomes
	\begin{eqnarray}\label{7p2}
		\begin{split}		
			&\int_0^T\int_M\varphi(1-f)\phi W^2\mu dt+\frac12\int_0^T\int_M[d\varphi(\nabla^{\nabla u}(\phi W))-\phi W\varphi_t]d\mu dt\\ 
			&+\int_0^T\int_M\varphi\left[\frac{f}{1-f}df(\nabla^{\nabla u}(\phi W))+\frac{1}{\phi}d\phi(\nabla^{\nabla u}(\phi W))\right]d\mu dt\\
			\leq&\int_0^T\int_M\varphi\left[fF_{\nabla u}(\nabla^{\nabla u}\phi)W^{\frac32}+\frac{F_{\nabla u}^2(\nabla^{\nabla u}\phi)}{\phi}W-\frac12(\Delta^{\nabla u}\phi-\phi_t)W\right.\\
			&\left.\quad\quad\quad\quad\quad\quad+K\phi W-\frac{q}{1-f}\phi W-\phi\frac{F_{\nabla u}(\nabla^{\nabla u}q)}{1-f}W^{\frac12}\right]d\mu dt,
		\end{split}
	\end{eqnarray}	
	To deal with the terms on the RHS of \eqref{7p2}, one could adopt the Young's inequality to obtain the following inequalities. That is,
	\begin{align}\label{7.21}
		\int_0^T\int_M\varphi fF_{\nabla u}(\nabla^{\nabla u}\phi)W^{\frac32}d\mu dt\leq& \frac16\int_0^T\int_M\varphi (1-f) \phi W^2d\mu dt\notag\\
		&+ C\int_0^T\int_M\varphi\frac{f^4}{(1-f)^3}\frac{F_{\nabla u}^4(\nabla \phi)}{\phi^3}d\mu dt,\\
		\int_0^T\int_M\varphi K\phi W d\mu dt\leq&\frac16\int_0^T\int_M\varphi (1-f)\phi W^2d\mu dt\notag\\
		&+C\int_0^T\int_M\varphi\frac{\phi K^2}{1-f}d\mu dt,\\
		\int_0^T\int_M\varphi \frac{q^-}{1-f}\phi Wd\mu dt\leq&\frac16\int_0^T\int_M\varphi (1-f)\phi W^2d\mu dt\notag\\
		&+C\int_0^T\int_M\varphi\frac{(q^-)^2}{(1-f)^3}\phi d\mu dt,\\
		\int_0^T\int_M\varphi \phi\frac{F_{\nabla u}(\nabla^{\nabla u}q)}{1-f}W^{\frac12}d\mu dt\leq&\frac16\int_0^T\int_M\varphi (1-f)\phi W^2d\mu dt\notag\\
		&+C\int_0^T\int_M\varphi\frac{F_{\nabla u}^{\frac43}(\nabla^{\nabla u}q)}{(1-f)^{\frac53}}d\mu dt,
	\end{align}
	and
	\begin{eqnarray}\label{7.25}
		\begin{split}
			\int_0^T\int_M\varphi\frac12\phi_t Wd\mu dt\leq&\frac16\int_0^T\int_M\varphi(1-f)\phi W^2d\mu dt\\
			&+C\int_0^T\int_M\varphi\frac{1}{1-f}\frac{\phi_t^2}{\phi}d\mu dt,
		\end{split}
	\end{eqnarray}
	where $C$ in these inequalities is a consistent constant. Plugging \eqref{7.21}-\eqref{7.25} and \eqref{est-J} into \eqref{7p2} yields	
	\begin{eqnarray}\label{7p3}
		\begin{split}		
			&\frac12\int_0^T\int_M[d\varphi(\nabla^{\nabla u}(\phi W))-\phi W\varphi_t]d\mu dt\\ 
			\leq&\int_0^T\int_M\varphi\frac{C_3}{R^2}[1+R(1+\sqrt{K(2R)+K'(2R)})]Wd\mu dt\\
			&-\int_0^T\int_M\varphi\left[\frac{f}{1-f}df(\nabla^{\nabla u}(\phi W))+\frac{1}{\phi}d\phi(\nabla^{\nabla u}(\phi W))\right]d\mu dt\\
			&-\frac16\int_0^T\int_M\varphi(1-f)\phi W^2\mu dt+C\int_0^T\int_M\varphi\left[\frac{f^4}{(1-f)^3}\frac{F_{\nabla u}^4(\nabla \phi)}{\phi^3}\right.\\
			&\left.\quad\quad\quad\quad\quad\quad+\frac{\phi K^2}{1-f}+\frac{(q^-)^2}{(1-f)^3}+\frac{F_{\nabla u}^{\frac43}(\nabla^{\nabla u}q)}{(1-f)^{\frac53}}+\frac{1}{1-f}\frac{\phi_t^2}{\phi}\right]d\mu dt\\
			=:&\int_0^T\int_M\varphi\Xi d\mu dt.
		\end{split}
	\end{eqnarray}	
	
	
	Suppose $(x_3,t_3)$ is the point so that $\phi W$ attains the maximum  in $Q_{2R,T}$. By the same argument as before, we know that $\Xi\geq0$ at $(x_3,t_3)$. Considering the peroperties iii) satisfied by the function $\tilde\phi$ and the fact that $f\leq 0$, we can deduced from \eqref{7p3} that
	\begin{eqnarray}\label{7p5}
		\begin{split}
			\phi^2W^2\leq& \frac{C_3}{R^2}[1+R(1+\sqrt{K(2R)+K'(2R)})](\phi W)\\
			&+C\left[\frac{1}{R^4}+\frac{1}{T^2}+K^2+\|q^-\|_{L^{\infty}(Q_{2R,T})}^2+\|F(\nabla q)\|_{L^{\infty}(Q_{2R,T})}^{\frac43}\right],
		\end{split}
	\end{eqnarray}
	where the constants $C_3$ and $C$ may be chosen larger than the ones in \eqref{7p3}. According to the inequality that $a_0\leq a_2+\sqrt{a_1}$, if $a_0^2\leq a_1+a_2a_0$, \eqref{7p5} implies that
	\begin{eqnarray}
		\begin{split}
			\phi W\leq& \frac{C_3}{R^2}[1+R(1+\sqrt{K(2R)+K'(2R)})]\\
			&+C\left[\frac{1}{R^2}+\frac{1}{T}+K+\|q^-\|_{L^{\infty}(Q_{2R,T})}+\|F(\nabla q)\|_{L^{\infty}(Q_{2R,T})}^{\frac23}\right].
		\end{split}
	\end{eqnarray}
\end{proof}

\section{Applications of the gradient estimates}\label{Ch-ap}

In this section, we will present some results about the solutions to the Finslerian Schr\"odinger equation \eqref{equ-FS1} by employing the gradient estimates. One may notice that when $q=0$, the Finslerian Schr\"odinger equation becomes the Finslerian heat equation. The Harnack inequality of the positive solution to the Finslerian heat equation and its mean value type inequality, as well as the bounds of the heat kernel have already been obtained by other researchers \cite{Ohta2014}\cite{XiaQ2020}\cite{Xia2023}.

\begin{thm}\label{thm-Harnack1-cpt}
	Let $(M,F,\mu)$ be a compact manifold without boundary, whose weighted Ricci curvature $Ric_N$ is bounded from below by $-K$, with $K\geq 0$. Let $q(x,t)$ be a function defined on $M\times[0,T]$ which is $C^2$ in the $x$-variable and $C^1$ in the $t$-variable. Assume that $\Delta^Vq\leq \theta$ and $F(\nabla q)\leq \gamma$ on $M\times [0,T]$ for some constant $\theta$ and $\gamma$. If $u(x,t)$ is a positive solution to the equation
	$$(\Delta^{\nabla u}-\partial_t-q)u(x,t)=0$$
	on $M\times (0,T]$, then for any $\beta>1$, $0<t_1<t_2\leq T$, and $x_1,x_2\in M$, we have the inequality
	\begin{eqnarray}
		u(x_1,t_1)\leq u(x_2,t_2)\left(\frac{t_2}{t_1}\right)^{\frac{N\beta}{2}}\exp (\mathcal{P}(t_2-t_1)+\mathcal{Q}_{\beta,R}(x_1,x_2,t_2-t_1)),
	\end{eqnarray}
	where
	\begin{eqnarray}
		\mathcal{P}=\mathcal{P}(K,\theta,\gamma,\beta,N)=C(N)\left[(\gamma^2(\beta-1)\beta)^{\frac13}+\frac{\beta}{\beta-1}K+(\beta\theta)^{1/2}\right],
	\end{eqnarray}
	and 
	\begin{eqnarray}
		\begin{split}
			&\mathcal{Q}_{\beta,R}(x_1,x_2,t_2-t_1)\\
			=&\inf_{\gamma\in\Gamma(M)}\left\{\frac{\beta}{4(t_2-t_1)}\int_0^1F^2(\dot\gamma)ds+(t_2-t_1)\int_0^1q(\gamma(s),st_1+(1-s)t_2)ds\right\},
		\end{split}
	\end{eqnarray}
	with the infimum taken over all forward paths parameterized by $[0,1]$ from $x_2$ to $x_1$, denoted by $\Gamma(M)$.
\end{thm}

The proof of Theorem \ref{thm-Harnack1-cpt} is the same to the one on noncompact forward complete Finsler metric measure spaces (Theorem \ref{thm-Harnack1-noncpt}). So we omit the proof here.  

\begin{thm}\label{thm-Harnack1-noncpt}
	Let $(M,F,\mu)$ be a complete manifold with boundary $\partial M$. Assume $p\in M$ and let $B_p(2R)$ be a forward geodesic ball of radius $2R$ centered at $p$ which does not intersect $\partial M$. We denote $-K(2R)$, with $K(2R)\geq0$, to be a lower bound of the mixed weighted Ricci curvatrue $^mRic^N$ on $B_p(2R)$ and denote $A(2R)\geq 1$ to be a local upper bound of the misalignment $\alpha$. We also suppose the non-Riemannian tensors satisfy $F(U)+F^*(\mathcal{T})+F(\mathrm{div} C(V))\leq K_0$,
	for any local vector fields $V,W$ in the definitions of $U$ and $\mathcal{T}$. Let $q(x,t)$ be a function defined on $M\times [0,T]$ which is $C^2$ in the $x$-variable and $C^1$ in the $t$-variable. Assume that $\Delta^Vq\leq \theta(2R)$ and $F(\nabla q)\leq \gamma(2R)$ on $B_p(2R)\times [0,T]$ for some constant $\theta(2R)$ and $\gamma(2R)$. If $u(x,t)$ is a positive solution to the equation
	$$(\Delta^{\nabla u}-\partial_t-q)u(x,t)=0$$
	on $M\times (0,T]$, then for any $\beta>1$, $0<t_1<t_2\leq T$, and $x_1,x_2\in B_p(R)$, we have the inequality
	\begin{eqnarray}\label{equ-thm-Harnack}
		u(x_1,t_1)\leq u(x_2,t_2)\left(\frac{t_2}{t_1}\right)^{\frac{N\beta}{2}}\exp (\mathcal{P}(t_2-t_1)+\mathcal{Q}_{\beta,R}(x_1,x_2,t_2-t_1)),
	\end{eqnarray}
	where
	\begin{eqnarray}
		\begin{split}
			\mathcal{P}=&\mathcal{P}(K,R,\theta,\gamma,\beta,A,N,K_0)=C(A,N,K_0)\\
			&\cdot\left[\frac{\beta}{R^{-1}}(\sqrt{K}+1)+\frac{\beta^3}{(\beta-1)R^{2}}+\left(\gamma^{2}(\beta-1)\beta\right)^{\frac{1}{3}}+(\beta\theta)^{1/2}+\frac{\beta}{\beta-1}K\right].
		\end{split}
	\end{eqnarray}
	and 
	\begin{eqnarray}
		\begin{split}
			&\mathcal{Q}_{\beta,R}(x_1,x_2,t_2-t_1)\\
			=&\inf_{\gamma\in\Gamma(R)}\left\{\frac{\beta}{4(t_2-t_1)}\int_0^1F^2(\dot\gamma)ds+(t_2-t_1)\int_0^1q(\gamma(s),st_1+(1-s)t_2)ds\right\},
		\end{split}
	\end{eqnarray}
	with the infimum taken over all forward paths in $B_p(R)$ parameterized by $[0,1]$ from $x_2$ to $x_1$, denoted by $\Gamma(R)$.
\end{thm}

\begin{proof}
	Let $\gamma$ be any curve by $\gamma:[0,1]\rightarrow B_p(R)$, with $\gamma(0)=x_2$ and $\gamma(1)=x_1$. We define $\eta:[0,1]\rightarrow B_p(R)\times [t_1,t_2]$ by
	$$\eta(s)=(\gamma(s),st_1+(1-s)t_2).$$
	Clearly, $\eta(0)=(x_2,t_2)$ and $\eta(1)=(x_1,t_1)$. Integrating $\frac{d}{ds}\log u$ along $\eta$, one gets that
	\begin{eqnarray}
		\begin{split}
			\log u(x_1,t_1)-\log u(x_2,t_2)&=\int_0^1\left(\frac{d}{ds}\log u\right)ds\\
			&=\int_0^1\{<\dot\gamma,\nabla\log u>-(t_2-t_1)(\log u)_t\}ds.
		\end{split}
	\end{eqnarray}
	Applying Theorem \ref{thm-GEcomplete-1} to $u$ yields
	\begin{eqnarray}\label{inequ-Hanack-1}
		\begin{split}
			&\log\left(\frac{u(x_1,t_1)}{u(x_2,t_2)}\right)\\
			\leq&\int_0^1\left\{F(\dot\gamma)F(\nabla\log u)-(t_2-t_1)\frac{F^2(\nabla \log u)}{\beta}+(t_2-t_1)\left[\mathcal{P}+\frac{N\beta}{2t}+q\right]\right\}ds\\
			\leq&\int_0^1\left\{\frac{\beta F^2(\dot\gamma)}{4(t_2-t_1)}+(t_2-t_1)\left[\mathcal{P}+\frac{N\beta}{2t}+q\right]\right\}ds.
		\end{split}
	\end{eqnarray}
	Plugging $t=st_1+(1-s)t_2$ into \eqref{inequ-Hanack-1} gives
	\begin{eqnarray}
		\begin{split}
			\log\left(\frac{u(x_1,t_1)}{u(x_2,t_2)}\right)\leq&\int_0^1\left\{\frac{\beta F^2(\dot\gamma)}{4(t_2-t_1)}+(t_2-t_1)q(\gamma(s),st_1+(1-s)t_2)\right\}ds\\
			&+\frac{N\beta}{2}\log\left(\frac{t_2}{t_1}\right)+\mathcal{P}(t_2-t_1).
		\end{split}
	\end{eqnarray}
	The result \eqref{equ-thm-Harnack} follows from taking exponential of the above inequality and then taking infimum.
\end{proof}

\noindent

Applying Theorem \ref{thm-GEcomplete-1-cor} instead in the above process on forward complete Finsler metric measure spaces, the same method provides that
\begin{thm}
	Let $(M,F,\mu)$ be a forward complete Finsler metric measure space without boundary. Suppose $u(x,t)$ is a positive solution on $M\times [0,T]$ to the Finslerian Schr\"odinger equation \eqref{equ-FS1}. Assume it admits finite misalignment $\alpha\leq\infty$, and the mixed weighted Ricci curvature $^mRic^N$ of $M$ is bounded from below by $-K$, with $K\geq 0$. Moreover, the non-Riemannian tensors satisfy $F(U)+F^*(\mathcal{T})+F(\mathrm{div} C(V))\leq K_0$,
	for any local vector fields $V,W$ in the definitions of $U$ and $\mathcal{T}$. Suppose there is a point $p\in M$, a constant $\theta$, and a function $\gamma(x,t)$, such that $F(\nabla q)\leq\gamma(r(x),t)$ and $\Delta^Vq\leq\theta$ on $M\times [0,T]$ with a constant $\theta$, for any locally nonvanishing vector field $V$, where $r(x)$ denotes the forward distance from $p$ to $x$. Then for any points $x,y\in M$, and $0<t_1<t_2\leq T$, the following estimates are valid:
	\begin{itemize}
		\item[(i)] If $K=0$ and $\lim_{r\rightarrow\infty}\frac{\gamma(r,t)}{r}\leq\iota$, for all $t\in[0,T]$, then
		\begin{eqnarray}
			\begin{split}
				u(x_1,t_1)\leq &u(x_2,t_2)\left(\frac{t_2}{t_1}\right)^{\frac{N}{2}}\\
				\cdot&\exp\left(C(N,\alpha,K_0)(\iota^{\frac23}+\theta^{\frac12})(t_2-t_1)+\mathcal{Q}(x_1,x_2,t_2-t_1)\right),
			\end{split}
		\end{eqnarray}
		where 
		\begin{eqnarray*}
			\begin{split}
				&Q(x_1,x_2,t_2-t_1)\\
				=&\inf_{\gamma\in\Gamma(M)}\left\{\frac{1}{4(t_2-t_1)}\int_0^1F^2(\dot\gamma)ds+(t_2-t_1)\int_0^1q(\gamma(s),st_1+(1-s)t_2)ds\right\}.
			\end{split}
		\end{eqnarray*}
	\item[(ii)] If $\gamma(r,t)\leq \gamma_0$ for some constant $\gamma_0$ in $M\times [0,T]$, then
	\begin{eqnarray}
		\begin{split}
			u(x_1,t_1)\leq& u(x_2,t_2)\left(\frac{t_2}{t_1}\right)^{\frac{N\beta}{2}}\\
			\cdot&\exp\left(C(N,\alpha,K_0)(\gamma_0^{\frac23}+\theta^{\frac12}+\frac{K}{\beta-1})(t_2-t_1)+\mathcal{Q}_{\beta}(x_1,x_2,t_2-t_1)\right),
		\end{split}
	\end{eqnarray}
	for all $\beta\in(1,2)$, where
	\begin{eqnarray*}
		\begin{split}
			&Q_{\beta}(x_1,x_2,t_2-t_1)\\
			=&\inf_{\gamma\in\Gamma(M)}\left\{\frac{\beta}{4(t_2-t_1)}\int_0^1F^2(\dot\gamma)ds+(t_2-t_1)\int_0^1q(\gamma(s),st_1+(1-s)t_2)ds\right\}.
		\end{split}
	\end{eqnarray*}
	with the infrimum taken over all forward paths in $M$ parameterized by $[0,1]$ joining $x_2$ to $x_1$.
		\end{itemize}
\end{thm}

For any equation on a manifold, an immediate corollary of the gradient estimate is the Liouville type theorem. The following theorem gives such an example about the elliptic Schr\"odinger equation. 

\begin{thm}
	Let $(M,F,\mu)$ be a forward complete Finsler metric measure space without boundary. Suppose it admits a finite misalignment $\alpha\leq\infty$, and the mixed weighted Ricci curvature $^mRic^N$ of $M$ is bounded from below by $-K$. Moreover, the non-Riemannian tensors satisfy $F(U)+F^*(\mathcal{T})+F(\mathrm{div} C(V))\leq K_0$,
	for any local vector fields $V,W$ in the definitions of $U$ and $\mathcal{T}$. Assume $q$ is a $C^2$ function defined on $M$ with $\Delta^V q\leq 0$ and $F(\nabla q)=o(r(x))$, where $r(x)$ denotes the forward distance from a fixed point $p$ to $x$. If $\inf q<0$, then the equation
	$$(\Delta^{\nabla u}-q)u(x)=0$$
	dose not admit a positive solution on $M$. 
\end{thm}

\begin{proof}
	Let $u(x)$ be a positive solution to $(\Delta^{\nabla u}-q)u=0$. Applying Theorem \ref{thm-GEcomplete-1-cor} to this time independent solution, we arrive at the estimate
	\begin{eqnarray}
		\frac{F^2(\nabla u)}{u^2}\leq \frac{N}{2t}+q.
	\end{eqnarray}
	Let $t\rightarrow\infty$, and evaluating at a point where $q<0$, we have a contradiction for $\inf q<0$.
\end{proof}

The Liouville type theorem of Theorem \ref{thm-GEcomplete-2} is about the ancient solution to \eqref{equ-FS1}. An ancient solution is a solution defined in all space and negative time.

\begin{thm}
	Let $(M,F,\mu)$ be a forward complete Finsler metric measure space without boundary. Assume it admits a finite misalignment $\alpha\leq\infty$, and the mixed weighted Ricci curvature $^mRic^{\infty}$ of $M$ is bounded from below by $-K$, with $K\geq 0$, the $S$-curvature has a positive $L^{\infty}$ upper bound $K'$. moreover, the non-Riemannian tensors satisfy $F(U)+F^*(\mathcal{T})+F(\mathrm{div} C(V))\leq K_0$,
	for any local vector fields $V,W$ in the definitions of $U$ and $\mathcal{T}$. Suppose $q(x,t)$ in \eqref{equ-FS1} satisfies that
	\begin{itemize}
		\item[i)] $q(x,t)\equiv q(x)$, i.e., $q$ is time independent;
		\item[ii)] $\|q^-\|_{L^{\infty}(B(p,R))}=o(R^{-2})$ when $R\rightarrow\infty$;
		\item[iii)] $\|\nabla q\|_{L^{\infty}(B(p,R))}=o(R^{-3})$ when $R\rightarrow\infty$.
	\end{itemize}
	Then
	\begin{itemize}
		\item[(1)] for $q(x)\not\equiv 0$, equation \eqref{equ-FS1} has no positive ancient solution with $\log u(x,t)=o(d(p,x)+\sqrt{\vert t\vert })$ near infinity;
		\item[(2)] for $q(x)\equiv 0$, equation \eqref{equ-FS1} becomes heat equation, and it has only constant positive ancient solution with $\log(u(x,t)+1)=o(d(p,x)+\sqrt{\vert t\vert })$ near infinity.
	\end{itemize}
\end{thm}

\begin{proof}
	To prove (1), suppose $u(x,t)$ is a positive ancient solution with $\log u(x,t)=o(d(p,x)+\sqrt{\vert t\vert })$ near infinity. Fixing $(x_0,t_0)$ in space-time and using Theorem \ref{thm-GEcomplete-2} for $u$ on $B(x_0, R)\times [t_0-R^2,t_0]$, we obtain that
	\begin{eqnarray*}
		\frac{F(\nabla u(x_0,t_0))}{u(x_0,t_0)}\leq \frac{C}{R}(1+o(R))
	\end{eqnarray*}
	near infinity. Letting $R\rightarrow \infty$ shows that $F(\nabla u(x_0,t_0))=0$. Such property holds for almost every $(x_0,t_0)$. We get $\nabla u\equiv 0$ for $u$ is in $H^1(M)$. Thus, $u(x,t)=u(t)$. Now we have arrived at that $q(x)$ is a constant, denoted by $\lambda$. The equation \eqref{equ-FS1} changes into 
	\begin{eqnarray}\label{du=lambda}
		\frac{d}{dt}u=\lambda u.
	\end{eqnarray}
	
	From ii) in the condition, we know that $\lambda\leq 0$. We will show that $\lambda$ must be 0 and get a contradiction with the hypothesis that $q(x)\not\equiv 0$. In fact, if $\lambda<0$, one could deduce by taking integral of \eqref{du=lambda} on $(t,0]$ with $t<0$ that
	\begin{eqnarray*}
		u(t)=u(0)e^{\lambda t},
	\end{eqnarray*}
	which is a contradiction with $u=e^{o(d(p,x)+\sqrt{\vert t\vert })}$ near infinity.
	
	To prove (2), let $u$ be the function satisfying $\log(u+1)=o(d(p,x),\sqrt{\vert t\vert })$ near infinity. Applying Theorem \ref{thm-GEcomplete-2} to $u+1$ on $B(x_0, R)\times [t_0-R^2,t_0]$ provides that 
	\begin{eqnarray*}
		\frac{F(\nabla u(x_0,t_0))}{u(x_0,t_0)+1}\leq \frac{C}{R}(1+o(R)).
	\end{eqnarray*}
	Letting $R\rightarrow 0$, it follows that $F(\nabla u(x_0,t_0))=0$. Since $(x_0,t_0)$ is arbitrary, one sees that $u$ is a constant.
\end{proof}

\section*{Acknowledge}
We are very grateful to the reviewers for their careful review and valuable comments. This work was partially supported by~NNSFC (No. 12001099, 12271093).


\hskip -0.6cm
Bin Shen\\
School of Mathematics, Southeast University, Nanjing 211189, P. R. China\\
E-mail: shenbin@seu.edu.cn

\end{document}